\documentclass[zamm,a4paper,fleqn]{w-art}
\usepackage{times}
\usepackage{w-thm}
\usepackage{cite}

\usepackage{graphicx,amssymb,amsmath,a4wide}
\def\dd{\,\mathrm{d}}
\def\expe{\mathrm{e}}
\newcommand{\real}{\mathbb{R}}
\newcommand{\nat}{\mathbb{N}}
\def\bpf{\begin{proof}}
\def\epf{\end{proof}}
\def\be{\begin{equation}\label}
\def\ee{\end{equation}}
\newtheorem{hypothesis}[theorem]{Hypothesis}
\newcommand{\eps}{\varepsilon}  
\renewcommand{\c}{\underline{c}}
\newcommand{\e}{e}
\newcommand{\calP}{{\mathcal{P}}}    
\newcommand{\calU}{{\mathcal{U}}}    
\newcommand{\calW}{{\mathcal{W}}}    
\newcommand{\calT}{{\mathcal{T}}}    

\usepackage[]{graphicx}
\begin{document}
\DOIsuffix{theDOIsuffix}
\Volume{VV}
\Issue{I}
\Month{MM}
\Year{YYYY}
\pagespan{3}{}
\Receiveddate{XXXX}
\Reviseddate{XXXX}
\Accepteddate{XXXX}
\Dateposted{XXXX}



\title[A thermodynamically consistent model for ferroelectricity and -elasticity]{A thermodynamically consistent phenomenological model for ferroelectric and ferroelastic hysteresis
\thanks{Supported by GA\v CR Grant P201/10/2315 and RVO: 67985840, as well as FWF grant P24970}
}

\author[B. Kaltenbacher]{Barbara Kaltenbacher\inst{1,}%
  \footnote{e-mail:~\textsf{barbara.kaltenbacher@aau.at} 
}
}
\address[\inst{1}]{Institute for Mathematics, Alpen-Adria-Universit\"at, Universit\"atsstr. 65-67, A-9020 Klagenfurt, Austria}
\author[P. Krej\v{c}\'{i}]{Pavel Krej\v{c}\'{i} \inst{2,} 
  \footnote{e-mail:~\textsf{krejci@math.cas.cz}}
}
\address[\inst{2}]{Institute of Mathematics, Czech Academy of Sciences,
\v{Z}itn\'{a} 25, CZ-11567 Praha 1, Czech Republic}
\begin{abstract}
We propose a hysteretic model for electromechanical coupling in piezoelectric materials, with the strain and
the electric field as inputs and the stress and the polarization as outputs. This constitutive
law satisfies the thermodynamic principles and exhibits good agreement with experimental measurements.
Moreover, when it is coupled with the mechanical and
electromagnetic balance equations, the resulting PDE system is well-posed under the hypothesis that hysteretic effects
take place only in one preferred direction. We prove the existence and uniqueness of its global weak solutions for each
initial data with prescribed regularity. One of the tools is a new Lipschitz continuity theorem for the inverse Preisach
operator with time dependent coefficients.
\end{abstract}
\maketitle                   






\section{Introduction}
The piezoelectric effect is a coupling between electrical and
mechanical fields within certain materials that has numerous
applications ranging from ultrasound generation in medical imaging and
therapy via acceleration sensors and injection valves in automotive
industry to high precision positioning systems. 
Driven by the increasing demand for devices operating at high field
intensities especially in actuator applications, the field of
hysteresis modeling for piezoelectric materials is currently one of
highly active research.  
The approaches that have been considered so far can be divided into basically four categories:
\begin{enumerate}
\item[(1)] {\em Thermodynamically consistent models}
being based on a macroscopic view to
describe microscopic phenomena in such a way
that the second law of thermodynamics is satisfied,
see for example, 
Bassiouny and Ghaleb \cite{bassiouny89:2}, 
Kamlah and B\"ohle \cite{kamlah01:1}, 
Landis \cite{landis04:1}, 
Schr\"oder and Romanowski \cite{Schroeder05}, 
Su and Landis \cite{Landis07}, 
Linnemann et al. \cite{Linnemann09}.

\item[(2)] {\em Micromechanical models}
that consider the material on the level of single grains, see, for example,
Delibas et al. \cite{Buelent05},  
Fr\"ohlich \cite{froehlich01:1}, 
Huber and Fleck \cite{Huber01},  
Belov and Kreher \cite{Belov06}, 
Huber \cite{Huber06}, 
McMeeking et al. \cite{McMeeking07},
Smith and Hu \cite{Smith12}.

\item[(3)] {\em Phase field models}
that describe the transition between phases (corresponding to the
motion of walls between domains with different polarization
orientation) using the Ginzburg-Landau equation  
for some order parameter, see, for example  
Xu et al \cite{Xuetal},
Wang et al \cite{Wangetal}.  

\item[(4)] {\em Phenomenological models using hysteresis operators}
partly originating from the input-output description of piezoelectric devices for control purposes,
see, for example, 
Hughes and Wen \cite{hughes95:1}, 
Kuhnen \cite{kuhnen01:1}, 
Cima et al. \cite{Cimaa02}, 
Smith et al. \cite{Smith03}, 
Ball et al. \cite{Ball07}, 
Pasco and Berry \cite{Pasco04}.
\end{enumerate}

Also multiscale coupling between macro- and microscopic as well as phase field
models partly even down to atomistic simulations has been investigated, see, e.g., 
\cite{MieheZaehRosato12,SchroederKeip10}.

For a mathematical analysis of some thermodynamically consistent models we refer
to \cite{AlberKraynyukova12,KraynyukovaNesenenko13,MielkeTimofte06}.

Whereas most of the so far existing models are designed for the simulation
of polarization, depolarization or cycling along the main hysteresis loop,
the simulation of actuators requires the accurate simulation of minor loops as well. 
Moreover, the physical behavior can so far be reproduced only qualitatively, whereas
the use of models in actuator simulation (possibly also aiming at simulation
based optimization) needs to match measurements precisely.
Simulation of a piezoelectric device with a possibly complex geometry requires not
only an input-output model but needs to resolve the spatial distribution of the
crucial electric and mechanical field quantities, which leads to partial
differential equations. Therewith, the question of numerical
efficiency becomes important.

Preisach operators \cite{BrokateSprekels,KrasPokr,KrejciBuch,mayergoyz91:2,VisintinBuch}
are phenomenological models for rate independent hysteresis that are capable of reproducing
minor loops and can be very well fitted to measurements. Moreover, they allow for a highly
efficient evaluation by the application of certain memory deletion rules and the use of
so-called Everett or shape functions.

Motivated by these facts, in \cite{HKKL08,KKHL10}, a model for ferroelectric hysteresis
under uniaxial loading using Preisach operators is proposed and studied.
Like most of the above mentioned models, it is based on an additive decomposition of
the strain and the  polarization into reversible and irreversible parts. The reversible
quantities follow the linear piezoelectric material law, while the irreversible polarization
is represented by a Preisach operator of the imposed electric field and the irreversible
strain is a polynomial function of the polarization. Moreover, the piezoelectric
coupling coefficient is proportional to the polarization.

However, for the model from \cite{HKKL08,KKHL10}, it is not clear whether or under which
conditions the second law of thermodynamics is satisfied. We will therefore here consider
a new  material law, which is inspired by the one proposed for magnetostriction \cite{DKV13}
whose thermodynamic consistency is based on the use of hysteresis potentials, and which
is additionally able to capture ferroelastic effects.

The paper is structured as follows. In Section \ref{mode}, we derive piezoelectric
constitutive equations involving a Preisach hysteresis operator from basic thermodynamic principles
and show some simulation results with the proposed model.
Section \ref{prei} is devoted to the proof of a new result on Lipschitz invertibility of the Preisach
operator with time dependent coefficients, and in Section \ref{long} we construct a unique solution
of a system of electromechanical balance equations by Banach contraction principle.
We restrict ourselves to the situation of uniaxial loading, hence scalar constitutive relations. 
A perspective to the situation of electric and mechanic fields depending on three space variables while still considering uniaxial loading, is provided in Section \ref{thin}.
The case of vector or actually tensor valued hysteresis will be subject of future research.

\section{The model}\label{mode}

Differently from most of the above cited approaches, we consider the electric field $E$
and the mechanical strain $\eps$ as state variables and the dielectric displacement
$D=D(\eps,E)$ as well as the mechanical stress $\sigma=\sigma(\eps,E)$ as state functions.
The reason for doing so lies in the fact that the balance equations

(namely Newton's law on the mechanical side and Gauss law on the electric side)
\begin{eqnarray}
\rho \ddot{u} -\nabla_s^T \sigma &=&0\label{Newton}\\
-\nabla\cdot D&=&0\label{Gauss}
\end{eqnarray}
(where $u$ is the mechanical displacement, $\nabla_s$ the symmetric gradient and $\nabla_s^T$ the dyadic divergence) are naturally formulated in terms of these state functions.

For the interrelation between these quantities,
we assume, similarly to \cite{DKV13}, that
hysteresis effects are due to one single Preisach operator $\calP$ with potential $\calU$
acting on an auxiliary state function $q = q(\eps,E)$.
In the scalar case, this leads to the assumption that the stress $\sigma$,
the dielectric displacement $D$, and the free energy $F=F(\eps,E)$ are of the form
\begin{eqnarray}
\sigma &=&\c \eps - \e E + a \calP[q] + b \calU[q],
\label{law_sigma_1d}\\ 
D&=& \e \eps +\kappa E + c \calP[q] + d \calU[q],
\label{law_D_1d}\\ 
F&=& \frac{\c}{2}\eps^2+\frac{\kappa}{2}E^2 + \xi \calP[q] + \eta \calU[q],
\label{ansatz_F_1d}
\end{eqnarray}
where the coefficients $a=a(\eps,E)$, $b=b(\eps,E)$, $c=c(\eps,E)$, $d=d(\eps,E)$,
$\xi=\xi(\eps,E)$, $\eta=\eta(\eps,E)$ as well as the function  $q=q(\eps,E)$,
are to be determined in agreement with the principles of thermodynamics.
The elastic, dielectric, and piezoelectric coupling coefficients 
$\c>0$, $\e\in\real$, and $\kappa>0$ are given constants.
It will be shown in Section \ref{prei} that the Preisach hysteresis operator
with nonnegative density satisfies the energy inequality
\begin{equation}\label{ccwhystpot}
q(t)\frac{\dd}{\dd t}\calP[q](t)-\frac{\dd}{\dd t}\calU[q](t)\geq 0\ \mbox{ a.e. }\quad 
\forall q\in W^{1,1}(0,T).
\end{equation}
Similarly to Section 3.4 in \cite{DKV13}, we now derive conditions to ensure
thermodynamic admissibility, namely the requirement that 
\begin{equation}\label{thermodyn_1d}
\dot{\eps}\sigma+\dot{D} E - \dot{F}\geq 0
\end{equation}
has to hold for all processes $\eps$, $\sigma$, $D$, $E$ that obey
the constitutive laws \eqref{law_sigma_1d}, \eqref{law_D_1d}. We obtain
\begin{equation}\label{ineq_law}
0\leq 
(a\dot{\eps}+E\dot{c}-\dot{\xi}) \calP[q]
+ (b\dot{\eps}+E\dot{d}-\dot{\eta}) \calU[q]
+ (Ec-\xi)\frac{\dd}{\dd t}\calP[q]\\ 
+ (Ed-\eta)\frac{\dd}{\dd t}\calU[q]\,.
\end{equation}
Both $\calP[q]$ and $\calU[q]$ are hysteresis operators and may take arbitrary values
independent of each other, so that we have to demand
\begin{equation}\label{eq_abcd_ccw}
\begin{aligned}
0 =&\ a\dot{\eps}+E\dot{c}-\dot{\xi} = \dot\eps\left(a+ E\frac{\partial c}{\partial \eps}
-\frac{\partial \xi}{\partial \eps}\right) + \dot E \left(E\frac{\partial c}{\partial E}
-\frac{\partial \xi}{\partial E}\right),\\
0=&\ b\dot{\eps}+E\dot{d} -\dot{\eta} = \dot\eps\left(b+ E\frac{\partial d}{\partial \eps}
-\frac{\partial \eta}{\partial \eps}\right) + \dot E \left(E\frac{\partial d}{\partial E}
-\frac{\partial \eta}{\partial E}\right).
\end{aligned}
\end{equation}
These relations have to hold for all processes, hence the coefficients of $\dot\eps$ and
$\dot E$ have to vanish identically. In other words, if we set $f_1=\eta-Ed$, $f_2 = Ec - \xi$,
we must have
\begin{equation}\label{coef}
a = - \frac{\partial f_2}{\partial \eps}, \ b = \frac{\partial f_1}{\partial \eps}, \
c = \frac{\partial f_2}{\partial E},\ d = -\frac{\partial f_1}{\partial E}.
\end{equation}
Inequality \eqref{ineq_law} then becomes
\begin{equation}\label{ineq_law_ccw}
0\leq f_2\frac{\dd}{\dd t}\calP[q]- f_1\frac{\dd}{\dd t}\calU[q]\,,
\end{equation}
and will be satisfied, by virtue of \eqref{ccwhystpot}, provided we choose $q := f_2/f_1$ with $f_1>0$.

We see that the whole model depends on the choice of two state functions
$f_1 = f_1(\eps,E)>0$ and $f_2 = f_2(\eps,E)$ which characterize the material properties.
A canonical choice, which is sufficient in many situations, consists in putting
$f_2(\eps,E) = E$, $f_1(\eps,E) = f(\eps)$ with a suitable function $f$ of one variable. Then the constitutive
law \eqref{law_sigma_1d}--\eqref{ansatz_F_1d} becomes
\begin{eqnarray}
\sigma &=&\c \eps - \e E +  f'(\eps) \calU\left[\frac{E}{f(\eps)}\right],
\label{law_sigma_1dr}\\ 
D&=& \e \eps +\kappa E + \calP\left[\frac{E}{f(\eps)}\right],
\label{law_D_1dr}\\ 
F&=& \frac{\c}{2}\eps^2+\frac{\kappa}{2}E^2 + f(\eps) \calU\left[\frac{E}{f(\eps)}\right].
\label{ansatz_F_1dr}
\end{eqnarray}
The full 1D system
for unknown functions $u(x,t)$, $E(x,t)$, $(x,t)\in (0,\ell)\times(0,T)$,
 describing longitudinal oscillations of a piezoelectric beam, then reads
\begin{equation}\label{PDEhyst}
\begin{aligned}
\rho u_{tt} - \sigma_x =&\ 0,\\
D_x=&\ 0,\\
\eps =&\ u_x.
\end{aligned}
\end{equation}
The equation $D_x = 0$ means that $D$ is a function of $t$ only, say, $D(x,t)=r(t)$, that is,
\begin{equation}\label{inverse}
\e \eps +\kappa E + \calP\left[\frac{E}{f(\eps)}\right] = r(t)\,,
\end{equation}
where $r(t)$ is a function which is known from the boundary condition
$D(0,t)=D(\ell,t)=r(t)$, corresponding to an impressed (or measured) boundary current.
Furthermore, we complement the mechanical constitutive law \eqref{law_sigma_1dr} with a viscosity term
$\nu \eps_t$, where $\nu>0$ is the viscosity coefficient, that is,
\begin{eqnarray}
\sigma &=&\nu \eps_t+ \c \eps - \e E +  f'(\eps) \calU\left[\frac{E}{f(\eps)}\right]
\label{law_sigma_1dv}
\end{eqnarray}
instead of \eqref{law_sigma_1dr}. We prescribe
some boundary conditions, for example $u=0$ on $x=0$, and $\sigma = s(t)$ on $x=\ell$,
which 
corresponds to the experimental setting of a beam which is clamped at the left tip, with an impressed (or measured) force at the right beam tip. 
We resume the analysis of this system
in Section \ref{long} and before, in Section \ref{prei}, we establish some new properties of the Preisach operator
which will enable us to eliminate $E$ from the system by solving Eq.~\eqref{inverse} independently
with respect to $q = E/f(\eps)$.

\begin{remark}
{\rm We wish to mention that independence of the dielectric saturation value
on the stress (which is the case e.g. for the setting $c=1$, $d=0$ considered
in Section 3.4.1 of \cite{DKV13}) physically makes sense also here: Note that
saturation of the polarization, i.e., of $P=\calP[q]$ taking its maximal absolute value,
corresponds to the situation of all c-axes (and therewith all elementary dipoles) of the
single crystals being aligned as much as possible to the load axis, under the given
geometric constraints (note that each grain has its own coordinate system of preferred
directions). This maximal alignment gives the same polarization, independent of the
imposed stress. In other words, no matter how large the imposed stress is, there exists
a sufficiently large impressed electric field (in load axis direction) that will bring
the polarization to its maximal value by aligning all elementary dipoles (and therewith
all c-axes) as much as possible in load direction.}
\end{remark}

\begin{figure}
\begin{center}
\includegraphics[width=0.4\textwidth]{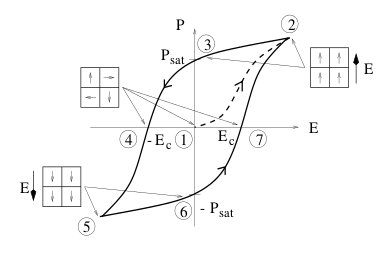}
\includegraphics[width=0.4\textwidth]{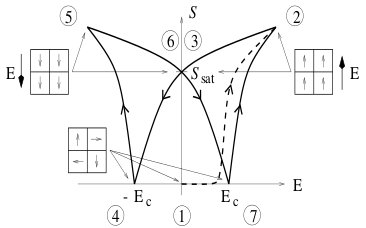}
\end{center}
\caption{polarization hysteresis (left) and strain (butterfly) hysteresis (right) under a bipolar electric field excitation; 
taken from \cite{kamlah01}
\label{fig:Kamlah1}}
\end{figure}

\begin{figure}
\begin{center}
\includegraphics[width=0.4\textwidth]{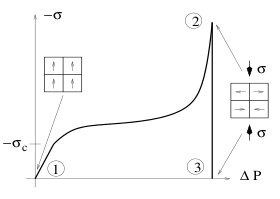}
\includegraphics[width=0.4\textwidth]{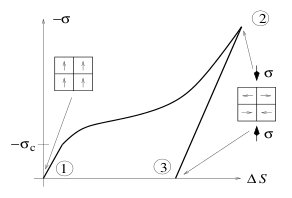}
\end{center}
\caption{mechanical depolarization (left) and stress-strain relation (right) under compressive stress load; 
taken from \cite{kamlah01}
\label{fig:Kamlah2}}
\end{figure}

\begin{figure}
\begin{center}
\includegraphics[width=0.6\textwidth]{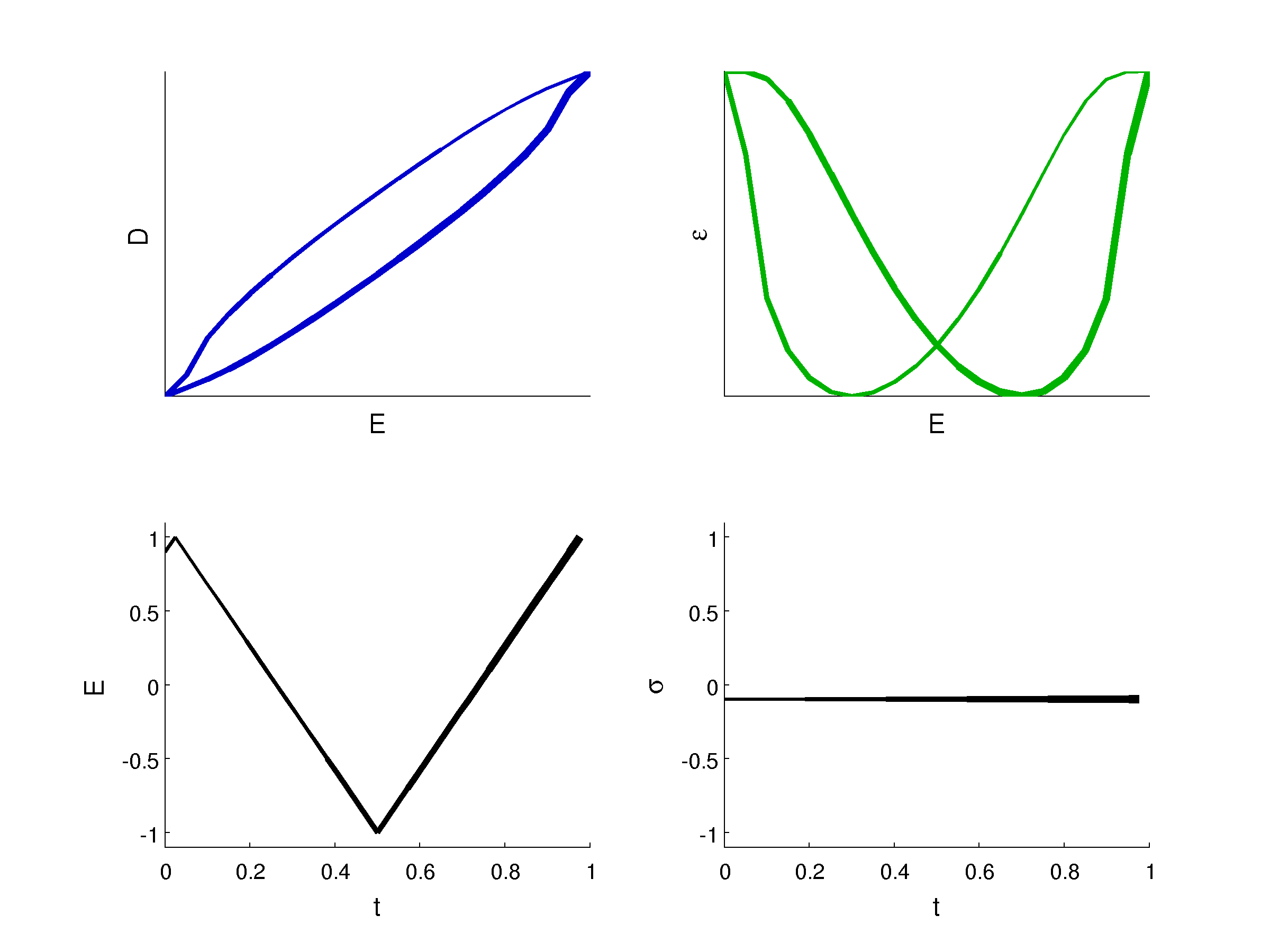}
\end{center}
\caption{polarization hysteresis (top row left) and strain (butterfly) hysteresis (top row right) under a bipolar electric field excitation (bottom row); 
simulations with model \eqref{ansatz_F_1dr} where $f$ is chosen according to \eqref{fcn_f1}, and $\calP$ according to Definition \ref{d1} with \eqref{fcn_g}.
\label{fig11}}
\end{figure}

\begin{figure}
\begin{center}
\includegraphics[width=0.6\textwidth]{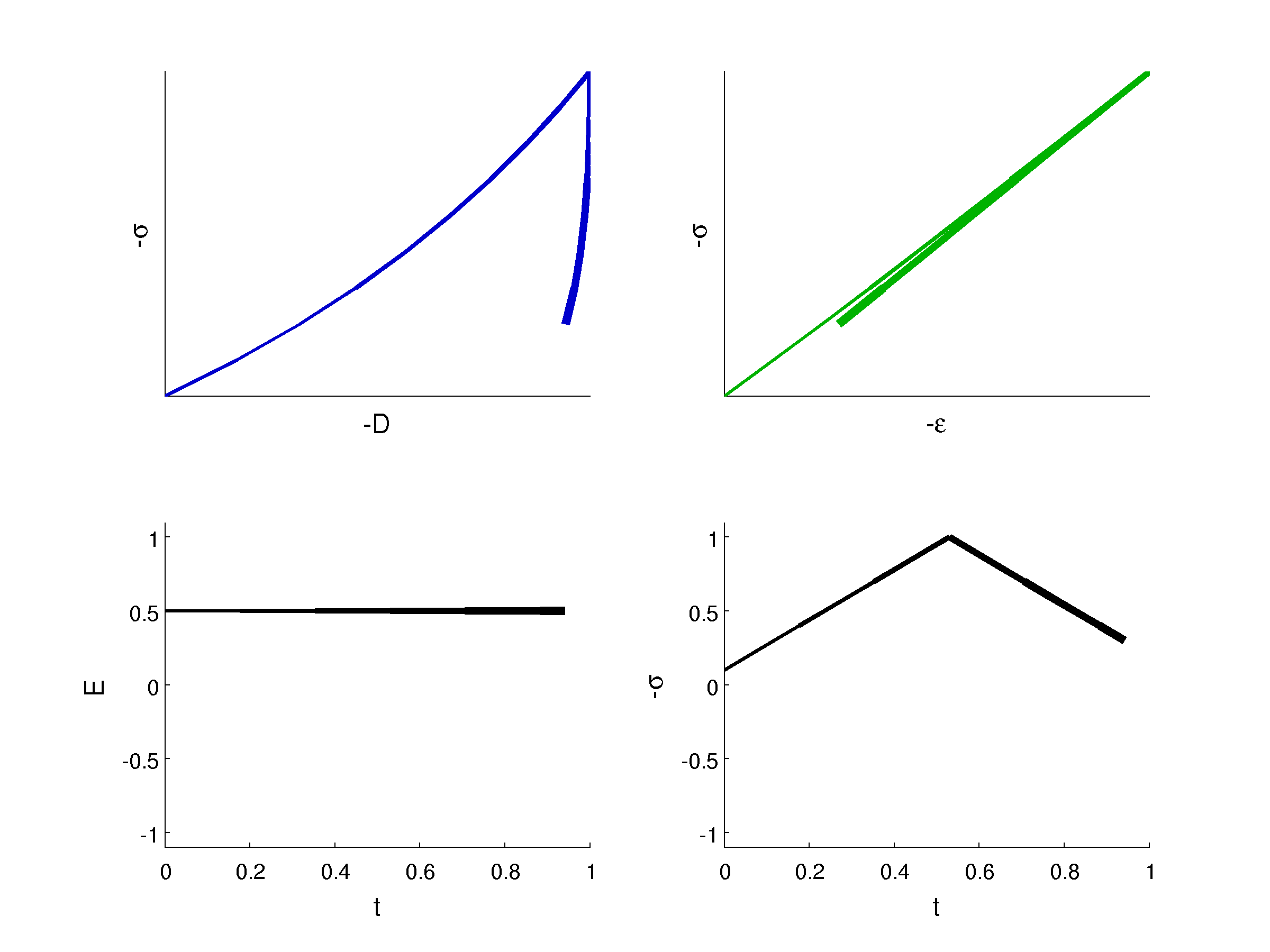}
\end{center}
\caption{mechanical depolarization (left) and stress-strain relation (right) under compressive stress load (bottom row); 
simulations with model \eqref{ansatz_F_1dr} where $f$ is chosen according to \eqref{fcn_f1}, and $\calP$ according to Definition \ref{d1} with \eqref{fcn_g}.
\label{fig21}}
\end{figure}

\begin{figure}
\begin{center}
\includegraphics[width=0.6\textwidth]{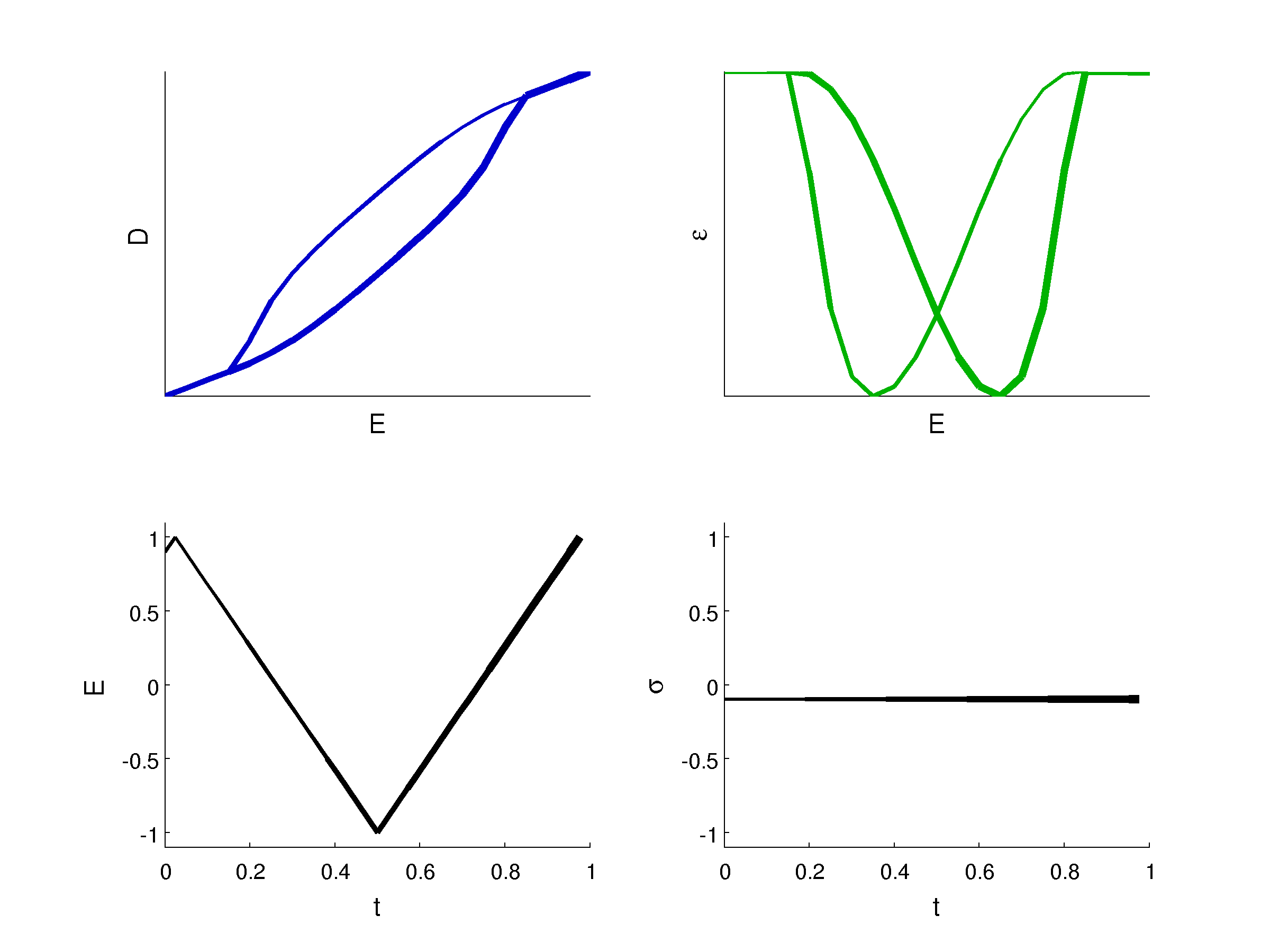}
\end{center}
\caption{polarization hysteresis (top row left) and strain (butterfly) hysteresis (top row right) under a bipolar electric field excitation (bottom row); 
simulations with model \eqref{ansatz_F_1dr} where $f$ is chosen according to \eqref{fcn_f4}, and $\calP$ according to Definition \ref{d1} with \eqref{fcn_g}.
\label{fig14}}
\end{figure}

\begin{figure}
\begin{center}
\includegraphics[width=0.6\textwidth]{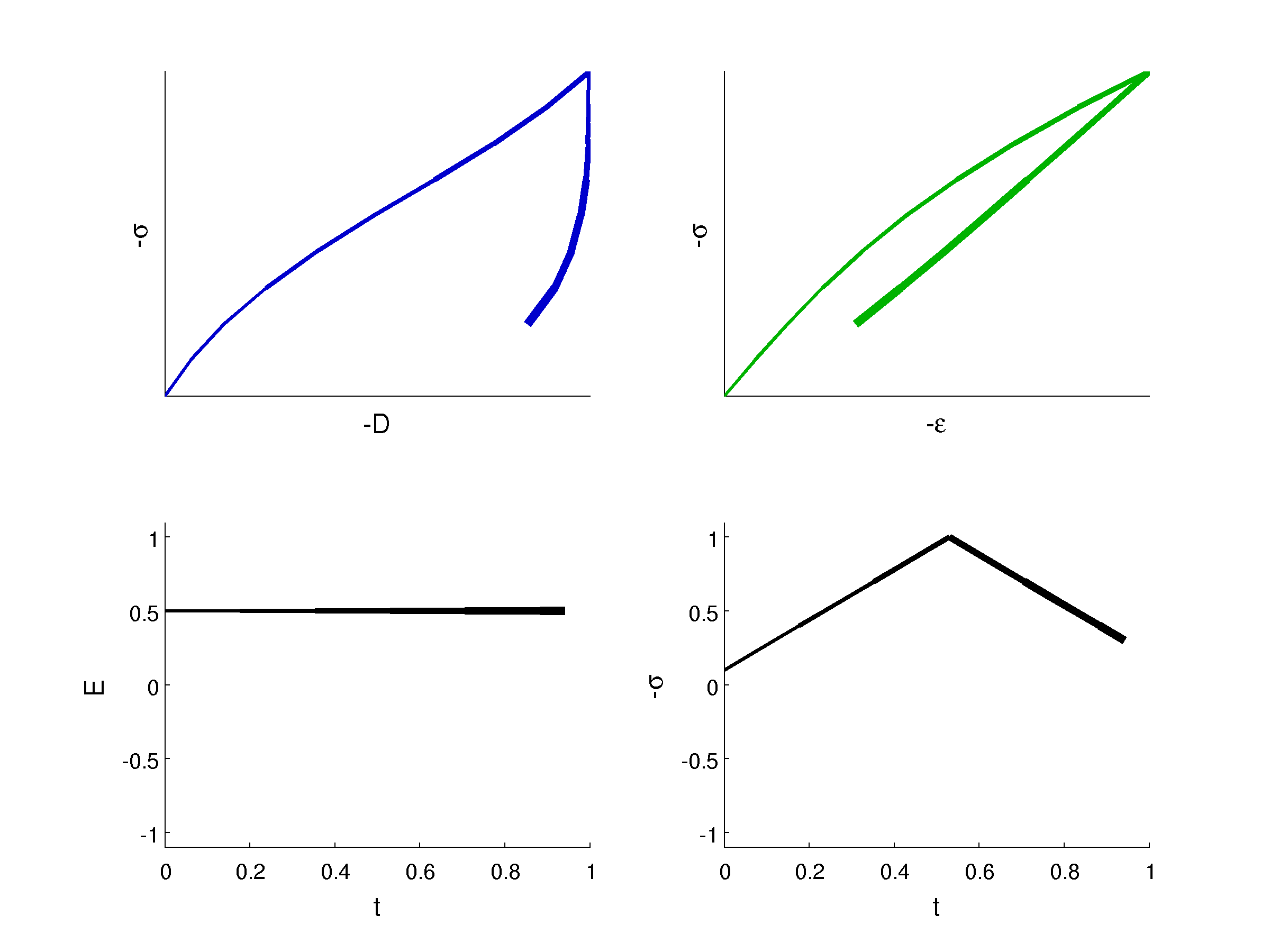}
\end{center}
\caption{mechanical depolarization (left) and stress-strain relation (right) under compressive stress load (bottom row); 
simulations with model \eqref{ansatz_F_1dr} where $f$ is chosen according to \eqref{fcn_f4}, and $\calP$ according to Definition \ref{d1} with \eqref{fcn_g}.
\label{fig24}}
\end{figure}

We conclude this section with some simulation results for a very simple choice of the Preisach operator $\calP$ and the function $f$, see Figures \ref{fig11}, \ref{fig21}, \ref{fig14}, \ref{fig24}, in order to illustrate that the expected qualitative behavior of ferroelectric and ferroelastic hysteresis (see Figures \ref{fig:Kamlah1}, \ref{fig:Kamlah2} taken from \cite{kamlah01}) can indeed be recovered. Here we use the functions  
\begin{equation}\label{fcn_g}
g(r,s)=\begin{cases}
\mbox{proj}_{[-1+r,1-r]}(s)= \max\{-1+r,\min\{1-r,s\}\}\mbox{ if }r\leq1\\
0 \mbox{ else }
\end{cases}
\end{equation}
in the Definition \ref{d1} of $\calP$, as well as 
\begin{equation}\label{fcn_f1}
f(x)=1.1-x\, \quad -1\leq x\leq 1
\end{equation}
or
\begin{equation}\label{fcn_f4}
f(x)=\frac12+\frac{(x-1)^4}{4}\, \quad -1\leq x\leq 1
\end{equation}
(in both cases the extension to $\real\setminus[-1,1]$ is done such that Hypothesis \ref{h1} is satisfied)
$\c=1$, $\e=0$, $\kappa=0.01$
in \eqref{ansatz_F_1dr},
and normalize all input quantities to the unit interval $[-1,1]$.
In Figures \ref{fig11}, \ref{fig21}, \ref{fig14}, \ref{fig24},increasing line thickness indicates proceeding time in order to show that the curves are traversed in the right direction.\\
Of course, also quantitative agreement with measurement can be expected to be achievable by approprate fitting methods.


\section{Inversion of hysteresis operators with time dependent coefficients}\label{prei}

In this section, we prove Lipschitz continuity statements for mappings which with given functions
$w, b \in C[0,T]$ associate the solution $q$ of the equation
\be{e1}
q(t) + b(t) \calP[q](t) = w(t) \quad \forall t\in [0,T]\,,
\ee
where $\calP: C[0,T] \to C[0,T]$ is a Lipschitz continuous operator in the sense that
\be{xe2}
|\calP[q_1](t) - \calP[q_2](t)| \le M \max_{\tau \in [0,t]}|q_1(\tau) - q_2(\tau)|
\ee  
for some $M>0$ and for all $q_1, q_2 \in C[0,T]$.

It follows from \eqref{xe2} that $\calP$ has the Volterra property, i.e., for all $t\in [0,T]$
\[ 
\forall q_1,q_2\in C[0,T]\ : \quad q_1=q_2 \mbox{ on } [0,t] \mbox{ implies }
\calP[q_1]=\calP[q_2] \mbox{ on } [0,t]\,.
\]
We first specify sufficient conditions under which is the mapping $(w,b) \mapsto q$ well defined.

\begin{lemma}\label{l1}
Let $\calP$ be as in \eqref{xe2}, and such that $I+c\calP$ admits for every $c\ge 0$ a Lipschitz continuous
inverse in the sense that the inequality
\be{xe4}
|(I+c\calP)^{-1}[w_1](t) - (I+c\calP)^{-1}[w_2](t)| \le L \max_{\tau \in [0,t]}|w_1(\tau) - w_2(\tau)|
\ee
holds for all $q_1, q_2 \in C[0,T]$ with Lipschitz constant $L>0$ independent of $c$.
Let $w,b \in C[0,T]$ be given functions, and let
$b(t) \ge 0$ for all $t\in [0,T]$. Then there exists a unique $q\in C[0,T]$ such that
\eqref{e1} holds.
\end{lemma}

The Lipschitz continuity of $(I + c\calP)^{-1}$ independently of the constant $c$ is fulfilled in many
typical situations arising in the theory of hysteresis operators. If $\calP$ is the Preisach operator
defined below in Definition \ref{d1}, it was shown in \cite{bv} that \eqref{xe4} is satisfied with $L=2$.
If $\calP$ belongs to the subclass of Prandtl-Ishlinskii operators, then $(I + c\calP)^{-1}$
admits in addition an explicit representation, see \cite{mathz}.

\bpf
We choose $\gamma< 1/(ML)$ and a division $0=t_0< t_1 <\dots < t_m=T$ of $[0,T]$
such that $|b(t) - b(t_{k-1})| < \gamma$ for $t\in [t_{k-1},t_k]$, $k=1, \dots, m$.
We prove by induction that
\[
\forall k\in\{1,\ldots,m\}\ \exists! q\in C[0,t_k]\ : \quad q\mbox{ satisfies \eqref{e1}}. 
\]

In $[0,t_1]$, we rewrite \eqref{e1} as
$$
q(t) + b(0) \calP[q](t) = w(t) - (b(t)-b(0))\calP[q](t)\,,
$$
or, equivalently,
\begin{equation}\label{e2}
q(t) = (I + b(0)\calP)^{-1}\big[w - (b(\cdot)-b(0))\calP[q]\big](t)\,.
\end{equation}
By our assumptions of Lipschitz continuity and the choice of $\gamma$,
the mapping on the right hand side of \eqref{e2} is a contraction
with respect to $q$, hence it admits a unique fixed point.

Assume now that we have constructed a unique $q(t)$ satisfying \eqref{e1} on $[0,t_{k-1}]$.
We introduce the set
$$
V_k = \{\hat q\in C[0,t_k]: \hat q(t) = q(t) \mbox{ on } [0,t_{k-1}]\}.
$$
We define a mapping $R_k: V_k \to V_k$ which with each $\hat q \in V_k$ associates
the solution $q$ of the equation
\begin{equation}\label{e3}
q(t) + b(t_{k-1}) \calP[q](t) = w(t) - (b(t)-b(t_{k-1}))\calP[\hat q](t)
\end{equation}
for $t \in [0,t_{k}]$. The mapping $R_k$ is again a contraction on $V_k$ 
(by our assumptions of Lipschitz continuity, the choice of $\gamma$, as well as the Volterra property)
and therefore admits a unique fixed point. The induction argument over $k$ then completes the proof.
\epf

Here, we focus on the case of Preisach operators. We use the following definition which is
shown in \cite{max} to be equivalent to the original Preisach construction in \cite{preisach}.

\begin{definition} \label{d1}
Let $g: (0,\infty) \times \real \to \real$ be a measurable function which is Lipschitz continuous
in the second variable and such that $g(r,0) = 0$ for a.\,e. $r>0$.
For a given input $q\in W^{1,1}(0,T)$,
we define the output $\calP[q] \in W^{1,1}(0,T)$ of a Preisach operator $\calP$ by the integral
\be{xe5}
\calP[q](t) = \int_0^\infty g(r, \xi_r(t))\dd r \quad \mbox{for } t \in [0,T]\,,
\ee
where $\xi_r \in W^{1,1}(0,T)$ is the unique solution of the variational inequality
\be{varin}
\left\{
\begin{array}{ll}
|q(t) - \xi_r(t)| \le r & \forall t\in [0,T]\,,\\
\dot\xi_r(t)(q(t) - \xi_r(t) - rz) \ge 0 & \mbox{a.\,e. } \ \forall |z|\le 1\,,\\
\xi_r(0) = \max\{q(0) - r, \min\{0, q(0) + r\}\}\,.&
\end{array}
\right.
\ee
The Preisach operator $\calP$ is called a Prandtl-Ishlinskii operator if $g$ is linear in $v$, that is,
$g(r,v) = \mu(r)v$ for some $\mu\in L^1_{loc}(0,\infty)$.
\end{definition}
We easily check from \eqref{varin} that $\xi_r(t) = 0$ for $r \ge \|q\|$, where
$\|\cdot \|$ denotes the sup-norm in $C[0,T]$, so that
the integral in \eqref{xe5} is meaningful.

The parameter $r$ is the memory variable, and the mapping $q \mapsto \xi_r$ introduced
in \cite{KrasPokr} is called the {\em play operator\/}. Note that its extension to $C[0,T]$ is
Lipschitz continuous with Lipschitz constant $1$, that is,
\be{xe5a}
|\xi_r^1(t) - \xi_r^2(t)| \le \max_{\tau \in [0,t]}|q_1(\tau) - q_2(\tau)|\,.
\ee

For our purposes, it is convenient to reduce the set of admissible functions $g$,
and we adopt the following hypothesis.

\begin{hypothesis}\label{h1}
We assume that $g(r,0) = 0$ for a.\,e. $r>0$, and
\begin{itemize}
\item[{\rm (i)}] $0 \le \frac{\partial g}{\partial v}(r,v) \le \mu(r)$ a.\,e., where
$\mu \in L^1(0,\infty)$, $\int_0^\infty \mu(r)\dd r =M$;
\item[{\rm (ii)}] There exists $\infty>M_1>0$ such that
\be{xe6}
\int_0^\infty\int_{-\infty}^{\infty} \frac{\partial g}{\partial v}(r,v)\dd v \dd r = M_1.
\ee
\end{itemize}
\end{hypothesis}
Hypothesis \ref{h1} (i) means that thr Preisach operator $\calP$ is dominated by a Prandtl-Ishlinskii operator.
Hypothesis \ref{h1} (ii) implies that $\|\calP[q]\|\leq M_1$ for all $q\in C[0,T]$.
It is easy to check that even without Hypothesis \ref{h1}\,(i), (ii), the Preisach operator satisfies the
energy inequality \eqref{ccwhystpot} with the choice
\be{xe7}
\calU[q](t) = \int_0^\infty G(r, \xi_r(t)) \dd r\,,
\ee
with
\be{xe8}
G(r,v) = \int_0^v v' \frac{\partial g}{\partial v}(r,v')\dd v'\,.
\ee
Moreover, it satisfies the Volterra property
by virtue of \eqref{xe5a}. Still by
\eqref{xe5a}, the operator $\calP$ can be extended to $C[0,T]$,
and if Hypothesis \ref{h1}\,(i) is fulfilled, then the Lipschitz property
\eqref{xe2} holds.

The main result of this section reads as follows.

\begin{theorem}\label{t1}
Let Hypothesis \ref{h1}\,(i) hold, and let $b, w_1, w_2 \in C[0,T]$ be given,
$0 \le b(t) \le \bar b$ for all $t\in [0,T]$. Let $q_1, q_2 \in C[0,T]$ be
such that
\be{xe9}
q_i(t) + b(t) \calP[q_i](t) = w_i(t) \quad \forall t\in [0,T]\,,\quad i=1,2\,.
\ee
Then we have
\be{xe9a}
\|q_1 - q_2\| \le \expe^{\bar b M}\|w_1 - w_2\|\,.
\ee
\end{theorem}

The proof of Theorem \ref{t1} is divided into several steps. Indeed, as $W^{1,1}(0,T)$
is dense in $C[0,T]$, it is enough to prove \eqref{xe9a} for $q_i \in W^{1,1}(0,T)$,
and we repeatedly use the variational inequality \eqref{varin}. The following observation
is due to Brokate, see \cite{BrokateSprekels}, and is related to the Brokate identity
in more general situations, cf.~\cite{adk}. We give a full proof here, also because
it uses techniques that will play a substantial role in the proof of Theorem \ref{t1}.

\begin{lemma}\label{l2}
Let $r_2 > r_1 > 0$ and $q \in W^{1,1}(0,T)$ be given, and let $\xi_{r_i}$, $i=1,2$, be the 
solution of \eqref{varin}. Then we have
\be{bi}
\left\{
\begin{array}{ll}
|\xi_{r_1}(t) - \xi_{r_2}(t)| \le r_2 - r_1 & \forall t\in [0,T]\,,\\
\dot\xi_{r_2}(t)(\xi_{r_1}(t) - \xi_{r_2}(t) - (r_2 - r_1)z) \ge 0 & \mbox{a.\,e. } \ \forall |z|\le 1\,.
\end{array}
\right.
\ee
\end{lemma}

\bpf
Let $\eta\in W^{1,1}(0,T)$ be the solution of the variational inequality
\be{xe10}
\left\{
\begin{array}{ll}
|\xi_{r_1}(t) - \eta(t)| \le r_2 - r_1 & \forall t\in [0,T]\,,\\
\dot\eta(t)(\xi_{r_1}(t) - \eta(t) - (r_2 - r_1)z) \ge 0 & \mbox{a.\,e. } \ \forall |z|\le 1\\
\eta(0) = \max\{\xi_{r_1}(0) - (r_2 - r_1), \min\{0, \xi_{r_1}(0) + (r_2 - r_1)\}\}\,.&
\end{array}
\right.
\ee
Putting in the above inequality 
$z = \frac{\xi_{r_1}(t\pm\delta) - \eta(t\pm\delta)}{r_2-r_1}$
 and letting
$\delta$ tend to $0+$ we obtain
\be{xe10a}
\dot\eta(t)(\dot\xi_{r_1}(t)-\dot\eta(t))
= 0 \quad \mbox{a.\,e.}
\ee
We have by \eqref{varin} for $r=r_1$ that $\dot\xi_{r_1}(t)(q(t) - \xi_{r_1}(t) - r_1 z) \ge 0$ a.\,e.
for all $|z|\le 1$, and \eqref{xe10a} implies the implication $\dot\eta(t)\ne 0 \Longrightarrow
\dot\eta(t)=\dot\xi_{r_1}(t)$. Hence,
\be{xe11}
\dot\eta(t)(q(t) - \xi_{r_1}(t) - r_1 z) \ge 0\ \mbox{a.\,e. }\  \forall |z|\le 1.
\ee
We have $|q(t) - \eta(t)| \le r_2$, and adding \eqref{xe11} to \eqref{xe10} we obtain
\be{xe13}
\dot\eta(t)(q(t) - \eta(t) - r_2 z) \ge 0\ \mbox{a.\,e. }\  \forall |z|\le 1.
\ee
By \eqref{varin} for $r=r_2$ that
\be{xe12}
\dot\xi_{r_2}(t)(q(t) - \xi_{r_2}(t) - r_2 z) \ge 0\ \mbox{a.\,e. }\  \forall |z|\le 1.
\ee
It is easy to check that $\eta(0) = \xi_{r_2}(0)$, hence, comparing \eqref{xe12} with \eqref{xe13},
we conclude that $\eta(t) = \xi_{r_2}(t)$ for all $t$, which we wanted to prove.
\epf

Next, we consider {\em discrete Preisach operators\/}
\be{xe15}
\calP_k[q] = \sum_{j=1}^{k} g_j(\xi_{r_j})\quad \mbox{for } k \in \nat\,,
\ee
corresponding to a sequence $0 = r_0 < r_1 < \dots$ and an associated sequence
$\{g_j\}_{j\in \nat}$ of nondecreasing functions.

\begin{proposition}\label{p2}
Let $b\in C[0,T]$ be as in Theorem \ref{t1}, and let $\{g_j\}_{j\in \nat}$ be a sequence of Lipschitz
continuous functions in $\real$ such that $g_j(0) = 0$ and $0 \le g_j'(v) \le \mu_j$ for all $j\in \nat$
and a.\,e. $v\in \real$, where $\mu_j\ge 0$ are constants.
Let $q_1, q_2, w_1, w_2 \in C[0,T]$ be such that
\be{xe16}
q_i(t) + b(t) \calP_k[q_i](t) = w_i(t) \quad \forall t\in [0,T]\,,\quad i=1,2
\ee
for some $k \in \nat$. Then
\be{indu}
\|q_1 - q_2\| \le \rho_k \|w_1 - w_2\|\,,
\ee
where $\rho_k = \Pi_{j=1}^k(1 + \bar b \mu_j)$.
\end{proposition}

\bpf
As mentioned above, the proof will be carried out for $q_1, q_2 \in W^{1,1}(0,T)$.
We proceed by induction over $k$. For $k=1$, Eq.~\eqref{xe16} reads
\be{xe17}
q_i(t) + b(t) g_1(\xi_{r_1}^i)(t) = w_i(t) \quad \forall t\in [0,T]\,,\quad i=1,2\,,
\ee
where $\xi_{r_1}^i$ is the solution of \eqref{varin} with input $q =q_i$ and memory level
$r=r_1$. In \eqref{varin} we choose $z = (q_2- \xi_{r_1}^2)/r_1$ for $i=1$ and
$z = (q_1- \xi_{r_1}^1)/r_1$ for $i=2$, and obtain
\begin{eqnarray}\label{xe18}
\dot\xi_{r_1}^1(t)(q_1(t) - \xi_{r_1}^1(t) - q_2(t) + \xi_{r_1}^2(t)) &\ge& 0\quad \mbox{a.\,e.}\,,\\
\label{xe19}
\dot\xi_{r_1}^2(t)(q_2(t) - \xi_{r_1}^2(t) - q_1(t) + \xi_{r_1}^1(t)) &\ge& 0\quad \mbox{a.\,e.}
\end{eqnarray}
The inequalities remain to hold when we multiply \eqref{xe18} by $g_1'(\xi_{r_1}^1(t))$ and 
\eqref{xe19} by $g_1'(\xi_{r_1}^2(t))$. We sum up the two and obtain
\be{xe20}
\frac{\dd}{\dd t} \left(g_1(\xi_{r_1}^1)- g_1(\xi_{r_1}^2)\right)(t)
\left((q_1 - q_2) - (\xi_{r_1}^1-\xi_{r_1}^2)\right)(t) \ge 0\quad \mbox{a.\,e.}
\ee
For $t\in [0,T]$ put
\be{pe8}
V_1(t) = \max\left\{\left(g_1(\xi_{r_1}^1)- g_1(\xi_{r_1}^2)\right)^2(t),
\left(\mu_1\|w_1 - w_2\|\right)^2\right\}.
\ee
We claim that
\be{xe21}
\dot V_1(t) \le 0 \quad \mbox{a.\,e.}
\ee
To prove \eqref{xe21}, we proceed by contradiction. Assume that there exists
a Lebesgue point $t \in (0,T)$ of both $\dot\xi_{r_1}^1, \dot\xi_{r_1}^2$
such that $\dot V_1(t) >0$. Then we have
\begin{eqnarray}\label{xe23}
&& |g_1(\xi_{r_1}^1) - g_1(\xi_{r_1}^2)|(t) > \mu_1 \|w_1 - w_2\|\,,\\[2mm]\label{xe24}
&& \frac{\dd}{\dd t}\left(g_1(\xi_{r_1}^1)- g_1(\xi_{r_1}^2)\right)(t)
\left(g_1(\xi_{r_1}^1)- g_1(\xi_{r_1}^2)\right)(t) >0\,.
\end{eqnarray}
It follows from \eqref{xe17}, \eqref{xe20}, and \eqref{xe24} that
\be{pe11}
\left(g_1(\xi_{r_1}^1) - g_1(\xi_{r_1}^2)\right)(t)
\left((w_1 - w_2)(t) - b(t)(g_1(\xi_{r_1}^1) - g_1(\xi_{r_1}^2))(t)
- (\xi_{r_1}^1 - \xi_{r_1}^2)(t)\right)\ge 0\,,
\ee
and since $g_1$ is monotone and
$\mu_1|\xi_{r_1}^1 - \xi_{r_1}^2| \ge |g_1(\xi_{r_1}^1) - g_1(\xi_{r_1}^2)|$, 
the estimate 
\[
\begin{aligned}
\left(g_1(\xi_{r_1}^1) - g_1(\xi_{r_1}^2)\right)(t)
(\xi_{r_1}^1 - \xi_{r_1}^2)(t)
=&\left|g_1(\xi_{r_1}^1) - g_1(\xi_{r_1}^2)\right|(t)\,
|\xi_{r_1}^1 - \xi_{r_1}^2|(t)\\
\geq& \frac{1}{\mu_1} \left(g_1(\xi_{r_1}^1) - g_1(\xi_{r_1}^2)\right)^2(t)
\end{aligned}
\]
holds, hence we obtain
\be{pe15}
|g_1(\xi_{r_1}^1) - g_1(\xi_{r_1}^2)|(t) \le \frac{\mu_1}{1+ b(t)\mu_1}
|w_1 - w_2|(t) \le \mu_1\|w_1 - w_2\|\,,
\ee
in contradiction with \eqref{xe23}. Hence, \eqref{xe21} holds, $V_1(t)$ is nonincreasing,
and we have by \eqref{pe8} for every $t \in [0,T]$ that
\be{xe25}
|g_1(\xi_{r_1}^1) - g_1(\xi_{r_1}^2)|(t)
\le \max\left\{\left|g_1(\xi_{r_1}^1)- g_1(\xi_{r_1}^2)\right|(0),
\mu_1\|w_1 - w_2\|\right\}.
\ee
The monotonicity of $g_1$ and of the initial value mapping in \eqref{varin} now yield
\[
|w_1-w_2|(0)=|q_1- q_2|(0)+b(0) \left|g_1(\xi_{r_1}^1)- g_1(\xi_{r_1}^2)\right|(0)\,,
\]
hence
$$
|g_1(\xi_{r_1}^1)- g_1(\xi_{r_1}^2)|(0) \le \mu_1|\xi_{r_1}^1- \xi_{r_1}^2|(0)
\le \mu_1|q_1- q_2|(0)
\le \mu_1|w_1- w_2|(0),
$$
and we conclude from \eqref{xe17} and \eqref{xe25} that
\be{xe26}
\|q_1 - q_2\| \le \bar b \|g_1(\xi_{r_1}^1) - g_1(\xi_{r_1}^2)\| + \|w_1 - w_2\|
\le (1+ \bar b \mu_1)\|w_1 - w_2\|,
\ee
and the first induction step for $k=1$ is done.

Let now $k>1$ be arbitrary and assume that \eqref{indu} holds for $k-1$ in place of $k$. Eq.~\eqref{xe16}
can be written as
\be{xe27}
q_i(t) + b(t) \calP_{k-1}[q_i](t) = w_i(t) - b(t)g_k(\xi_{r_k}^i(t))
\quad \forall t\in [0,T]\,,\quad i=1,2\,,
\ee
and by induction hypothesis we have
\be{xe28}
\|q_1 - q_2\| \le  \rho_{k-1}\left(\|w_1 - w_2\|+ \bar b \|g_k(\xi_{r_k}^1) - g_k(\xi_{r_k}^2)\|\right).
\ee
By the choices 
$z=\frac{q_i(t)-\xi_{r_k}^i(t)}{r_k}$, $i=1,2$ in \eqref{varin}
and 
$z=\frac{\xi_{r_j}^i(t)-\xi_{r_k}^i(t)}{r_k-r_j}$, $i=1,2$ in Lemma \ref{l2}, we obtain similarly as in \eqref{xe18}--\eqref{xe19}
\begin{eqnarray}\label{xe18a}
\dot\xi_{r_k}^1(t)\left((q_1(t) - q_2(t)) - (\xi_{r_k}^1(t)-\xi_{r_k}^2(t))\right) &\ge& 0\quad \mbox{a.\,e.}\,,\\
\label{xe19a}
\dot\xi_{r_k}^2(t)\left((q_2(t) - q_1(t)) - (\xi_{r_k}^2(t)-\xi_{r_k}^1(t))\right) &\ge& 0\quad \mbox{a.\,e.}\,,\\
\label{xe18b}
\dot\xi_{r_k}^1(t)\left((\xi_{r_j}^1(t)-\xi_{r_j}^2(t)) - (\xi_{r_k}^1(t)-\xi_{r_k}^2(t))\right)
&\ge& 0\quad \mbox{a.\,e.}\ \ \mbox{for } j=1, \dots, k-1\,,\\
\label{xe19b}
\dot\xi_{r_k}^2(t)\left((\xi_{r_j}^2(t)-\xi_{r_j}^1(t)) - (\xi_{r_k}^2(t)-\xi_{r_k}^1(t))\right)
&\ge& 0\quad \mbox{a.\,e.}\ \ \mbox{for } j=1, \dots, k-1\,.
\end{eqnarray}
The same argument as in the transition from \eqref{xe18}--\eqref{xe19} to \eqref{xe20} yields
\begin{eqnarray}\label{xe20a}
\frac{\dd}{\dd t} \left(g_k(\xi_{r_k}^1)- g_k(\xi_{r_k}^2)\right)(t)
\left((q_1 - q_2) - (\xi_{r_k}^1-\xi_{r_k}^2)\right)(t) &\ge& 0\quad \mbox{a.\,e.},\\[2mm]\label{xe20b}
\frac{\dd}{\dd t} \left(g_k(\xi_{r_k}^1)- g_k(\xi_{r_k}^2)\right)(t)
\left((\xi_{r_j}^1-\xi_{r_j}^2) - (\xi_{r_k}^1-\xi_{r_k}^2)\right)(t) &\ge& 0\quad \mbox{a.\,e.}
\end{eqnarray}
for every $j=1, \dots, k-1$ in \eqref{xe20b}. We continue as above and define the function
\be{pe8a}
V_k(t) = \max\left\{\left(g_k(\xi_{r_k}^1)- g_k(\xi_{r_k}^2)\right)^2(t),
\left(\mu_k\|w_1 - w_2\|\right)^2\right\}
\ee
with the goal to prove that
\be{xe21a}
\dot V_k(t) \le 0 \quad \mbox{a.\,e.}
\ee
Assume that there exists
a Lebesgue point $t \in (0,T)$ of both $\dot\xi_{r_k}^1, \dot\xi_{r_k}^2$
such that $\dot V_k(t) >0$. Then we have
\begin{eqnarray}\label{xe23a}
&& |g_k(\xi_{r_k}^1) - g_k(\xi_{r_k}^2)|(t) > \mu_k \|w_1 - w_2\|\,,\\[2mm]\label{xe24a}
&& \frac{\dd}{\dd t}\left(g_k(\xi_{r_k}^1)- g_k(\xi_{r_k}^2)\right)(t)
\left(g_k(\xi_{r_k}^1)- g_k(\xi_{r_k}^2)\right)(t) >0\,.
\end{eqnarray}
It follows from \eqref{xe20a}--\eqref{xe20b} that
\begin{eqnarray}\label{xe20c}
\left(g_k(\xi_{r_k}^1)- g_k(\xi_{r_k}^2)\right)(t)
\left((q_1 - q_2) - (\xi_{r_k}^1-\xi_{r_k}^2)\right)(t) &\ge& 0,\\[2mm]\label{xe20d}
\left(g_k(\xi_{r_k}^1)- g_k(\xi_{r_k}^2)\right)(t)
\left((\xi_{r_j}^1-\xi_{r_j}^2) - (\xi_{r_k}^1-\xi_{r_k}^2)\right)(t) &\ge& 0\,, \ j=1, \dots k-1\,.
\end{eqnarray}
As a consequence of \eqref{xe20d} and the monotonicity of $g_k$, we have
\be{xe20e}
\big(g_k(\xi_{r_k}^1)- g_k(\xi_{r_k}^2)\big)(t)
\big(\xi_{r_j}^1-\xi_{r_j}^2\big)(t) \ge 0 \quad \mbox{for } j=1, \dots k-1\,,
\ee
and the monotonicity of $g_j$ yields in turn
\be{xe20f}
\big(g_k(\xi_{r_k}^1)- g_k(\xi_{r_k}^2)\big)(t)
\big(g_j(\xi_{r_j}^1)-g_j(\xi_{r_j}^2)\big)(t) \ge 0 \quad \mbox{for } j=1, \dots k-1\,.
\ee
Combining \eqref{xe20c} with \eqref{xe20f} we obtain
\be{xe30}
\big(g_k(\xi_{r_k}^1)- g_k(\xi_{r_k}^2)\big)(t)
\left(\left(q_1 + b\sum_{j=1}^{k-1}g_j(\xi_{r_j}^1)\right)
-\left(q_2 + b\sum_{j=1}^{k-1}g_j(\xi_{r_j}^2)\right) - (\xi_{r_k}^1- \xi_{r_k}^2)\right)(t) \ge 0 \,,
\ee
which is nothing but (cf. \eqref{xe27})
\be{pe11a}
\left(g_k(\xi_{r_k}^1) - g_k(\xi_{r_k}^2)\right)(t)
\left((w_1 - w_2)(t) - b(t)(g_k(\xi_{r_k}^1) - g_k(\xi_{r_k}^2))(t)
- (\xi_{r_k}^1 - \xi_{r_k}^2)(t)\right)\ge 0\,.
\ee
Using the monotonicity of $g_k$ and inequality $\mu_k|\xi_{r_k}^1 - \xi_{r_k}^2| \ge |g_k(\xi_{r_k}^1) - g_k(\xi_{r_k}^2)|$,
we thus have like in \eqref{pe15}
\be{pe15b}
|g_k(\xi_{r_k}^1) - g_k(\xi_{r_k}^2)|(t) \le \frac{\mu_k}{1+ b(t)\mu_k}
|w_1 - w_2|(t) \le \mu_k\|w_1 - w_2\|
\ee
in contradiction with \eqref{xe23a}. Hence, \eqref{xe21a} holds, $V_k(t)$ is nonincreasing,
and we have by \eqref{pe8a} for every $t \in [0,T]$ that
\be{xe25a}
|g_k(\xi_{r_k}^1) - g_k(\xi_{r_k}^2)|(t)
\le \max\left\{\left|g_k(\xi_{r_k}^1)- g_k(\xi_{r_k}^2)\right|(0),
\mu_k\|w_1 - w_2\|\right\}.
\ee
The monotonicity of $g_k$ and of the initial value mapping in \eqref{varin} now yield
$$
|g_k(\xi_{r_k}^1)- g_k(\xi_{r_k}^2)|(0) \le \mu_k|\xi_{r_k}^1- \xi_{r_k}^2|(0)
\le \mu_k|q_1- q_2|(0)\le \mu_k|w_1- w_2|(0),
$$
and we conclude from \eqref{xe28} and \eqref{xe25a} that
\be{xe26a}
\|q_1 - q_2\| \le \rho_{k-1}(1+ \bar b \mu_k)\|w_1 - w_2\| = \rho_{k}\|w_1 - w_2\|,
\ee
which completes the induction argument.
\epf

We are now ready to prove Theorem \ref{t1}.

\bpf{ (Theorem \ref{t1})}\
Let $b, q_1, q_2, w_1, w_2$ be as in \eqref{xe9}. We choose $R > \max\{\|q_1\|, \|q_2\|\}$, so that
$\xi_r^i(t) = 0$ for all $r\ge R$, $t\in [0,T]$, and $i=1,2$.
For every $\delta > 0$, we choose a partition $0 = r_0 < r_1 < \dots < r_m = R$ with
$r_j - r_{j-1} < \delta$.
Our goal is to approximate $\calP$ by $\calP_m$ of the form \eqref{xe15} by putting
$$
g_j(v) = \int_{r_{j-1}}^{r_j} g(r,v)\dd r\quad \mbox{for } v \in \real, \ j=1, \dots, m\,.
$$
We have for $i=1,2$ that
$$
|\calP[q_i] - \calP_m[q_i]| = \left|\int_0^\infty g(r, \xi_r^i)\dd r - \sum_{j=1}^m g_j(\xi_{r_j}^i)\right|
= \left|\sum_{j=1}^m\int_{r_{j-1}}^{r_j} (g(r, \xi_r^i) -  g(r,\xi_{r_j}^i))\dd r\right|.
$$
It follows from Lemma \ref{l2} that $|g(r, \xi_r^i) -  g(r,\xi_{r_j}^i)| \le \mu(r)(r_j - r_{j-1})$
for $r \in [r_{j-1}, r_j]$, hence
$$
|\calP[q_i] - \calP_m[q_i]|(t) \le \delta \int_0^R \mu(r) \dd r \le M\delta\,.
$$
We rewrite \eqref{xe9} as
\be{xe9b}
q_i(t) + b(t) \calP_m[q_i](t) = w_i(t) + (\calP_m[q_i]- \calP[q_i])(t) \quad \forall t\in [0,T]\,,
\quad i=1,2\,.
\ee
From Proposition \ref{p2} it follows that
$$
\|q_1 - q_2\| \le \rho_m\left(\|w_1 - w_2\| + 2M\delta\right)
$$
with $\rho_m=\Pi_{j=1}^m(1 + \bar b \mu_j)$, 
$\mu_j=\int_{r_{j-1}}^{r_j} \mu(r) \dd r$. 
Note that $\log \rho_m = \sum_{j=1}^m \log(1+\bar b \mu_j) \le \sum_{j=1}^m \bar b\mu_j \le \bar b M$.
Hence,
$$
\|q_1 - q_2\| \le \expe^{\bar b M}\left(\|w_1 - w_2\| + 2M\delta\right).
$$
We can choose $\delta$ arbitrarily small, and the assertion follows.
\epf

We now prove the Lipschitz continuous dependence of the inverse mapping
on $b$ under suitable assumptions.

\begin{proposition}\label{p3}
Let Hypotheses \ref{h1}\,(i), (ii) hold. Let $b_1, b_2  \in C[0,T]$

be such that $0\le b_i(t) \le \bar b$ for all $t\in [0,T]$, $i=1,2$,
and let $w_1, w_2 \in C[0,T]$ be given. Let 
$q_1, q_2 \in C[0,T]$ be solutions of the equations
\begin{equation}\label{e8}
q_i(t) + b_i(t) \calP[q_i](t) = w_i(t)\ \ \forall t\in [0,T]\,, \ i=1,2\,.
\end{equation}
Then we have
\begin{equation}\label{e9}
\|q_1 - q_2\| \le \expe^{\bar b M}(\|w_1 - w_2\| + M_1\|b_1 - b_2\|)\,.
\end{equation}
\end{proposition}

\bpf
We have 
\begin{equation}\label{e10}
\begin{aligned}
q_1(t) + b_1(t) \calP[q_1](t) =&\ w_1(t) \,,\\
q_2(t) + b_1(t) \calP[q_2](t) =&\ w_2(t) - (b_2(t)-b_1(t))\calP[q_2](t)\,,
\end{aligned}
\end{equation}
and Theorem \ref{t1} together with \eqref{xe6} yields
$$
\|q_1 - q_2\| \le \expe^{\bar b M}(\|w_1 - w_2\| + \|\calP[q_2]\|\,\|b_1 - b_2\|)
\le \expe^{\bar b M}(\|w_1 - w_2\| + M_1\|b_1 - b_2\|),
$$
which we wanted to prove.
\epf


\section{Longitudinal oscillations of a piezoelectric beam}\label{long}

We keep Hypothesis \ref{h1} on the Preisach operator $\calP$, and assume moreover
that the mapping $(r,v) \mapsto v \frac{\partial g}{\partial v}(r,v)$
belongs to $L^1((0,\infty)\times \real)$, and there exists $\mu_1 \in L^1(0,\infty)$ such that
$$
|v|\frac{\partial g}{\partial v}(r,v) \le \mu_1(r) \quad \mbox{a.\,e.}
$$
This guarantees that the potential operator $\calU$ of the form \eqref{xe7}
is bounded and Lipschitz continuous. The constitutive function $f(\eps)$
in \eqref{law_sigma_1dr}--\eqref{ansatz_F_1dr} will be assumed to possess the following properties.
 
\begin{hypothesis}\label{h4}
The function $f: \real \to \real$ has Lipschitz continuous derivative $f'$ and the functions $\eps\mapsto 1/f(\eps)$ are $\eps\mapsto\eps/f(\eps)$ are Lipschitz continuous as well.
\end{hypothesis}
Under these hypotheses, we reformulate the constitutive equation \eqref{law_sigma_1dv}
in the form
\begin{equation}\label{e11}
\sigma = \nu\eps_t 
+\c\eps
+\calW[\eps]
\end{equation}
with a Lipschitz continuous operator $\calW: C[0,T] \to C[0,T]$. Indeed,
Eq.~\eqref{inverse} can be rewritten as
\begin{equation}\label{e12}
q + \frac{1}{\kappa f(\eps)} \calP[q] = \frac{1}{\kappa f(\eps)}(r - \e \eps)\,, \ q = \frac{E}{f(\eps)}.
\end{equation}
The functions $\eps\mapsto 1/f(\eps), \eps\mapsto \eps/f(\eps)$ are Lipschitz continuous by Hypothesis \ref{h4}. Hence, by Proposition \ref{p3}, the mapping $\eps\mapsto q$ is Lipschitz continuous, and \eqref{law_sigma_1dv}
can be written as
\begin{equation}\label{e13}
\sigma = \nu\eps_t + \c \eps - \frac{\e}{\kappa}(r - \e \eps - \calP[q]) +  f'(\eps) \calU[q]\,,
\end{equation}
so that \eqref{e11} holds with
$$
\calW[\eps] = 
- \frac{\e}{\kappa}(r - \e \eps - \calP[q]) +  f'(\eps) \calU[q]\,.
$$
This enables us to state the PDE problem \eqref{PDEhyst} in the form
\begin{equation}\label{e14}
\rho u_{tt} 
- \nu u_{xxt} 
- \c u_{xx}
= \calW[u_x]_x
\end{equation}
and we couple it with boundary conditions
\begin{equation}\label{bc}
u(0,t) = 0\,, \ (\nu u_{xt} 
+\c u_{x}
+\calW[u_x])(\ell, t) = s(t)\,,
\end{equation}
and initial conditions
\begin{equation}\label{ini}
u(x,0) = u^0(x)\,, \ u_t(x,0) = u^1(x)\,.
\end{equation}
In variational form, the problem reads
\begin{equation}\label{vari}
\rho \int_0^\ell u_{tt}\phi \dd x + \int_0^\ell(
\c u_{x}+
\nu u_{xt} + \calW[u_x])\phi_x\dd x
= s(t)\phi(\ell)\quad a.e. \ \forall \phi \in X\,,
\end{equation}
where we set $X = \{\phi \in  W^{1,2}(0,\ell): \phi(0) = 0\}$.

\begin{theorem}\label{main}
Let $r \in C[0,T]$, $s\in L^2(0,T)$, $u^0\in X$, and $u^1\in L^2(0,\ell)$ be given, Then Problem
\eqref{ini}--\eqref{vari} admits a unique solution $u\in C([0,T];X)$ such that
$u_{xt} \in L^2((0,\ell)\times(0,T))$, 
$u_{t} \in C([0,T]; L^2(\Omega))$,
$u_{tt} \in L^2(0,T; X')$.
\end{theorem}
\bpf
For $v\in C([0,T];X)$ such that $v_{xt} \in L^2((0,\ell)\times(0,T))$ and $v(x,0) = u^0(x)$,
we find $u$ with the desired regularity as the solution of the linear problem
\begin{equation}\label{varia}
\rho \int_0^\ell u_{tt}\phi \dd x + \int_0^\ell(
\c u_{x}+
\nu u_{xt} + \calW[v_x])\phi_x\dd x
= s(t)\phi(\ell)\quad a.e. \ \forall \phi \in X
\end{equation}
with initial conditions \eqref{ini}. We now prove that the mapping $v \mapsto u$ is a contraction
in the space
$$
Y = \{v\in C([0,T];X): v_{xt} \in L^2((0,\ell)\times(0,T)), v(x,0) = u^0(x), v(0,t) = 0\}\,,
$$
endowed with a suitable norm defined below in \eqref{test4}.
Let $v,\hat v \in Y$ be given, and let $u, \hat u$ be the corresponding
solutions. We test the difference of Eqs.~\eqref{varia} for $u$ and $\hat u$ by $\phi = u_t - \hat u_t$
and obtain, using the Lipschitz continuity of the 
the linear part containing $\c$ and of
$\calW$, that
\begin{equation}\label{test}
\frac{\rho}{2}\frac{\dd}{\dd t} \int_0^\ell (u_{t}- \hat u_t)^2\dd x + \nu\int_0^\ell(u_{xt}- \hat u_{xt})^2
\dd x \le C\int_0^\ell |u_{xt}- \hat u_{xt}| \max_{\tau\in [0,t]}|v_x(x,\tau) - \hat v_x(x,\tau)|\dd x
\end{equation}
with a constant $C>0$. We have $ \max_{\tau\in [0,t]}|v_x(x,\tau) - \hat v_x(x,\tau)|
\le \int_0^t |v_{xt}- \hat v_{xt}|\dd \tau$, so that \eqref{test} can be further estimated
using H\"older's inequality as
\begin{equation}\label{test1}
\tilde\rho\frac{\dd}{\dd t} \int_0^\ell (u_{t}- \hat u_t)^2\dd x + \int_0^\ell(u_{xt}- \hat u_{xt})^2 \dd x
\le Ct\int_0^t \int_0^\ell (v_{xt}- \hat v_{xt})^2 (x,\tau)\dd x \dd\tau
\end{equation}
with some constant $\tilde\rho>0$ and a possibly larger constant $C>0$.
This is an inequality of the form
\begin{equation}\label{test2}
\dot\alpha(t) + \beta(t) \le Ct\int_0^t \delta(\tau)\dd\tau
\end{equation}
with
$$
\alpha = \tilde\rho\int_0^\ell (u_{t}- \hat u_t)^2\dd x\,, \ \beta = \int_0^\ell(u_{xt}- \hat u_{xt})^2 \dd x\,, \
\delta = \int_0^\ell (v_{xt}- \hat v_{xt})^2 \dd x\,.
$$
We multiply \eqref{test2} by $\expe^{-Ct^2}$ and obtain
$$
\frac{\dd}{\dd t}\left(\expe^{-Ct^2}\left(\alpha(t) + \frac12\int_0^t \delta(\tau)\dd \tau\right)\right)
+ 2Ct\expe^{-Ct^2}\alpha(t) + \expe^{-Ct^2}\beta(t) \le \frac12 \expe^{-Ct^2}\delta(t)\,.
$$
We now integrate the above inequality from $0$ to $T$ and conclude that the mapping $v \mapsto u$
is a contraction in $Y$ endowed with norm
\begin{equation}\label{test4}
\|v\| = \int_0^T \int_0^\ell \expe^{-Ct^2} |v_{xt}(x,t)|^2 \dd x\dd t\,,
\end{equation}
which implies existence and uniqueness of solutions.
\epf

\section{Thin structures under uniaxial loading}\label{thin}

To some extent, the results from Section \ref{long} can be extended to a spatially three dimensional setting, as long as the setup still allows to assume that the third component of the dielectric displacement is constant in polarization direction. 
Namely we first of all lift the model \eqref{law_sigma_1dr}--\eqref{ansatz_F_1dr} to the spatially 3-d setting by prescribing a fixed polarization direction $p$, onto which we project the electric field.
For simplicity we assume that 
$p$ is constant.
Now $\sigma$, $\eps$, $D$, $E$ are tensor and vector valued functions, respectively, and $\c$, $\e$ $\kappa$ are constant 4th, 3rd, and 2nd order material tensors ($\c$ and $\kappa$ are symmetric positive definite), but the internal variable $q$ is a scalar valued function of the strain $\eps$ and of the projected electric field $p\cdot E$. 
Moreover, $\calP$ is still a scalar Preisach hysteresis operator with counterclockwise hysteresis potential $\calU$, i.e., we assume \eqref{ccwhystpot} to hold. 
For some strictly positive scalar valued function $f:\real^6\to\real^+$ satisfying Hypothesis \ref{h4} (with the preimage space $\real$ replaced by $\real^6$), we consider the analog of \eqref{law_sigma_1dr}--\eqref{ansatz_F_1dr}
\begin{eqnarray}
\sigma &=&\c \eps - \e E + \calU[q] Df(\eps)
\label{law_sigma_r}\\
D&=&  \e^T \eps +\kappa E + \calP[q] p
\label{law_D_r}\\
F&=& \frac{1}{2}\eps:(\c\eps)+\frac{1}{2}E\cdot(\kappa E) + f(\eps) \calU[q]
\label{ansatz_F_r}\\
q&=&\frac{p\cdot E}{f(\eps)}\,,\nonumber
\end{eqnarray}
where $Df(\eps)$ is the gradient of $f$, $Df:\real^6\to\real^6$.
It is readily checked that this ensures thermodynamic admissibility, i.e., the higher dimensional analog of \eqref{thermodyn_1d}.

We now consider Eqs.~\eqref{Newton}--\eqref{Gauss} in the form
\begin{eqnarray}
\rho u_{tt} - \nabla_s^T \sigma &=&0 \label{Newton_r}\\
-\nabla\cdot D &=&0 \label{Gauss_r}\\
\eps &=&\nabla_s u \label{u}\\
E &=&\nabla\phi\,, \label{phi}
\end{eqnarray}
where $u$ is the mechanical displacement, $\phi$ the negative of the electric potential, and $\nabla_s$ the symmetric gradient.

Without loss of generality (and actually consistently with the usual notation) we assume that the polarization direction is parallel to the $z$ axis, i.e., 
$$p=e_z\,,$$ 
and that the tensor of dielectric coefficients takes the form $\kappa=\left(\begin{array}{cc}\kappa^{xy}&0\\0&\kappa^z\end{array}\right)$ with some positive definite $2\times2$ matrix $\kappa^{xy}$ and $\kappa^z>0$.
\\
Our main restriction is the assumption that
the $z$ component of the dielectric displacement does not change in $z$ direction
\begin{equation}\label{Dz0}
\frac{\partial}{\partial z} D^z=0\,,
\end{equation}
that the domain takes the 
cylindrical form 
\[
\Omega=\Omega^{xy}\times (0,\ell)
\]
and that the boundary conditions 
\begin{eqnarray}
&&D\cdot n=0 \mbox{ on }\{z\}\times\partial \Omega^{xy} \quad \forall z\in(0,\ell) \label{bcDxy}\\
&&D\cdot n=-D^z=-r \mbox{ on }\{0\}\times\Omega^{xy}, \quad
D\cdot n=D^z=r \mbox{ on } \{\ell\}\times\Omega^{xy}
\label{bcDz}
\end{eqnarray}
hold, i.e., a current is prescribed, whose average over the boundary vanishes.
(Note that therewith the electric potential can only be expected to be unique up to a constant, but the electric field will still be unique.)
Assumption \eqref{Dz0} is realistic, e.g., when dealing with a thin structure extended in the $xy$ plane and excited by imposing some prescribed normal current via a pairs of opposite electrodes on top and bottom. 
Thus we have $D^z(x,y,z,t)=r(x,y,t)$ for all $(x,y,z)\in\Omega$, $t\in(0,T)$.
Therewith, the combination of \eqref{law_D_r} with \eqref{Gauss_r} can be split as follows:
\begin{eqnarray}
\nabla^{xy}\cdot \Bigl((\e^T \eps)^{xy} +\kappa^{xy} E^{xy}\Bigr)&=&0
\label{Dxy}\\
(\e^T \eps)^{z} +\kappa^{z} E^{z}+\calP[q]&=&r
\label{Dz}
\end{eqnarray}
where the latter equation can be cast in a form convenient for application of the results from Section \ref{prei}
$$ 
q + \frac{1}{\kappa^z f(\eps)} \calP[q] = \frac{r-(\e^T \eps)^{z}}{\kappa^z f(\eps)}
$$
with $q=\frac{E^z}{f(\eps)}$.
Hence we can write 
$$
E^z(x,y,z,\cdot)=\Phi^z[\eps(x,y,z,\cdot)]
\quad \forall (x,y,z)\in\Omega
$$
with a Lipschitz continuous mapping 
$$
\Phi^z:C[0,T]\to C[0,T]\,.
$$
On the other hand, for the $xy$ part, testing \eqref{Dxy} with $\phi$ (so that $\nabla^{xy} \phi= E^{xy}$) integrating by parts with respect to $(x,y)$, and using \eqref{bcDxy}, we obtain the estimate 
\[
\|E^{xy}(\cdot,\cdot,z,t)\|_{L^2(\Omega^{xy})}\leq 
\frac{|\e|}{\lambda_{\min}(\kappa^{xy})}
\|\eps(\cdot,\cdot,z,t)\|_{L^2(\Omega^{xy})}
\quad \forall z\in(0,\ell), t\in(0,T)\,,
\]
where $\lambda_{\min}(\kappa^{xy})$ is the smallest eigenvalue of $\kappa^{xy}$. Thus, by integrating the square of both sides with respect to $z$
\[
\|E^{xy}(\cdot,\cdot,\cdot,t)\|_{L^2(\Omega)}\leq 
\frac{|\e|}{\lambda_{\min}(\kappa^{xy})}
\|\eps(\cdot,\cdot,\cdot,t)\|_{L^2(\Omega)}
\quad \forall t\in(0,T)\,.
\]
By linearity, this provides us with Lipschitz continuity (with constant $L_{\Phi^{xy}}=\frac{|\e|}{\lambda_{\min}(\kappa^{xy})}$) of the mapping 
$$
\Phi^{xy}:L^2(\Omega)^6\to L^2(\Omega)^2
$$
such that
$$
E^{xy}(\cdot,\cdot,\cdot,t)=\Phi^{xy}(\eps(\cdot,\cdot,\cdot,t))
\quad \forall t\in(0,T)\,.
$$
Analogously to the proof of Theorem \ref{main}, under Hypothesis \ref{h4}, we will rewrite \eqref{law_sigma_r}
in the form 
\begin{equation}\label{f11}
\sigma = \nu\eps_t + \c \eps + \calW(\eps)
\end{equation}
with the operator $\calW$ mapping the (tensor valued) function $\eps$ of $(x,y,z,t)$ to a (tensor valued) function $\calW(\eps)$ of $(x,y,z,t)$ as follows:
\begin{align*}
\calW(\eps)(x,y,z,t)= &
- e \left(\begin{array}{c}
\Phi^{xy}(\eps(\cdot,\cdot,\cdot,t))(x,y,z)\\
\Phi^{z}[\eps(x,y,z,\cdot)](t)
\end{array}\right) \\
&+
\calU\left[\frac{\Phi^{z}[\eps(x,y,z,\cdot)]}{f(\eps(x,y,z,\cdot))} \right]
Df(\eps(x,y,z,t))
\end{align*}
To further proceed along the lines of the proof of Theorem \ref{main}, after elimination of the electric field we can now consider a purely mechanical problem
\begin{eqnarray}
\rho u_{tt}-\nabla_s^T\Bigl(\c\nabla_s u+\nu\nabla_s u_t+\calW(\nabla_s u)\Bigr)=0 \ \mbox{ in }\Omega\times(0,T)\\
u(x,y,z,0)=u^0\,, \quad u_t(x,y,z,0)=u^1 \ \mbox{ for all }(x,y,z)\in \Omega \label{ini3d}\\
u=0\ \mbox{ on }\Gamma\,, \quad 
n\cdot\Bigl(\nu\eps_t+\calW(\nabla_s u)\Bigr)=s \ \mbox{ on }\partial\Omega\setminus\Gamma\,, \label{bcu}
\end{eqnarray}
i.e., clamped at the boundary part $\Gamma$  and traction free or loaded with some surface force on the remainder of the boundary.

\bigskip

The test space $X$ becomes
$X = \{\phi \in  W^{1,2}(\Omega)^3: \phi = 0\ \mbox{ on }\Gamma\}$, the ansatz space 
$$
Y = \{v\in C([0,T];X): \nabla_s v_t \in L^2(\Omega\times(0,T))^6, \
v(t=0) = u^0, \ v = 0\ \mbox{ on }\Gamma\}\,,
$$
and the fixed point mapping $\calT$ maps $v\in Y$ to the unique solution $u\in Y$ of the linear variational problem
\begin{equation}\label{varia3d}
\rho \int_\Omega u_{tt}\phi \dd (x,y,z) + \int_\Omega(\c\nabla_s u
+\nu \nabla_s u_t + \calW(\nabla_s v)): \nabla_s \phi\dd (x,y,z)
= \int_\Gamma s\phi\dd S \quad a.e. \ \forall \phi \in X
\end{equation}
with initial conditions \eqref{ini3d}, provided $s\in L^2(0,T;H^{-1/2}(\Gamma))$. 
Again, to show contractivity, for $v,\hat{v}\in Y$, we test the difference between the equations for $u=\calT v$ and $\hat{u}=\calT \hat{v}$ with $u_t-\hat{u}_t$ and use Young's inequality to end up with the estimate
\begin{align*}
&\frac{\rho}{2}\frac{\dd}{\dd t} \|u_t- \hat u_t\|_{L^2(\Omega)}^2
+\frac{1}{2}\frac{\dd}{\dd t} \|\sqrt{\c}(\nabla_s u- \nabla_s \hat u)\|_{L^2(\Omega)}^2
+ \frac{\nu}{2}\|\nabla_s u_t- \nabla_s \hat u_t\|_{L^2(\Omega)}^2\\
&\le \frac{1}{2\nu} \|\calW(\nabla_s v)-\calW(\nabla_s \hat{v})\|_{L^2(\Omega)}^2\,.
\end{align*}
It remains to show that $\calW$ obeys the Lipschitz condition
\[
\max_{\tau\in[0,t]}\|\calW(\eps)-\calW(\hat{\eps})\|_{L^2(\Omega)}^2
\le C t \int_0^t \|\eps_t-\hat{\eps}_t\|_{L^2(\Omega)}^2 \dd \tau
\]
for some constant $C>0$ and any $\eps,\hat{\eps}\in L^2(\Omega\times(0,T))^6$ such that $\eps_t,\hat{\eps}_t\in L^2(\Omega\times(0,T))^6$ and $\eps(t=0)=\hat{\eps}(t=0)$. For the terms containing $\Phi^z$ and $\calU$, this follows analogously to the 1-d case from the respective Lipschitz continuity. For the term containing $\Phi^{xy}$ we have 
\begin{align*}
&\max_{\tau\in[0,t]}\int_\Omega |\Phi^{xy}(\eps(\cdot,\cdot,\cdot,\tau))(x,y,z)-
\Phi^{xy}(\hat{\eps}(\cdot,\cdot,\cdot,\tau))(x,y,z)|^2 \dd (x,y,z)\\
&\le L_{\Phi^{xy}}^2 \max_{\tau\in[0,t]}\|\eps-\hat{\eps}\|_{L^2(\Omega)}^2
= L_{\Phi^{xy}}^2 \max_{\tau\in[0,t]}\|\int_0^\tau(\eps_t-\hat{\eps}_t)\dd \tau_1\|_{L^2(\Omega)}^2
\le L_{\Phi^{xy}}^2 \, t \int_0^t\|\eps_t-\hat{\eps}_t\|_{L^2(\Omega)}^2 \dd \tau\,.\end{align*}
The rest of the proof goes exactly like in the 1-d case.

Therewith we have
\begin{corollary} \label{cor}
Let $r \in C([0,T];L^2(\Omega^{x,y}(0)\cup\Omega^{x,y}(\ell)))$,
$s\in L^2(0,T;H^{-1/2}(\Gamma))$, $u^0\in X$, and $u^1\in L^2(\Omega)$ be given. Then Problem
\eqref{law_sigma_r}--\eqref{phi} with boundary conditions \eqref{Dz0}--\eqref{bcDz}, \eqref{bcu} and initial conditions \eqref{ini3d}
admits a unique solution $(E^{xy},E^z,u)\in 
C([0,T];L^2(\Omega)) \times L^2(\Omega;C[0,T])\times C([0,T];X)$ such that 
$\nabla_s u_t \in L^2(\Omega\times(0,T))^6$, 
$u_{t} \in C([0,T]; L^2(\Omega)^3)$,
$u_{tt} \in L^2(0,T; X')$,
where $X = \{\phi \in  W^{1,2}(\Omega)^3: \phi = 0\ \mbox{ on }\Gamma\}$.
\end{corollary}

\bpf
Starting from the existence, uniqueness, and regularity of $u$, which we have already obtained along the lines of the proof of Theorem \ref{main}, we use Lipschitz continuity of the mappings $\Phi^{xy}$ and $\Phi^z$ to obtain assertions on $E^{xy}$, $E^z$. In the latter case, the estimate for obtaining the claimed regularity looks as follows:
\begin{align*}
&\int_\Omega \max_{t\in[0,T]} E^z(x,y,z,t)|^2 \dd (x,y,z)
\le \int_\Omega L_{\Phi^z} \max_{t\in[0,T]} |\nabla_s u(x,y,z,t)|^2 \dd (x,y,z)\\
&\le \int_\Omega L_{\Phi^z} \left| u^0(x,y,z)+\int_0^T \nabla_s u_t(x,y,z,t)\dd t\right|^2 \dd (x,y,z)\\
&\le 2 L_{\Phi^z} \Bigl( \|u^0\|_{L^2(\Omega)}^2+T \|\nabla_s u_t\|_{L^2(\Omega\times(0,T))}^2\Bigr)
\end{align*}
\epf

\begin{acknowledgement}
Support by GA\v CR Grant P201/10/2315 and RVO: 67985840, as well as FWF grant P24970 is gratefully acknowledged.
\end{acknowledgement}

%
\providecommand{\WileyBibTextsc}{}
\let\textsc\WileyBibTextsc
\providecommand{\othercit}{}
\providecommand{\jr}[1]{#1}
\providecommand{\etal}{~et~al.}


\begin{thebibliography}{[10]}

\bibitem{AlberKraynyukova12}
 \textsc{H.\,D. Alber} and  \textsc{N.~Kraynyukova},
A doubly nonlinear problem associated with a mathematical model for
  piezoelectric material behavior,
 \jr{ZAMM} \textbf{92}, 141--159 (2012).


\bibitem{Ball07}
 \textsc{B.\,L. Ball},  \textsc{R.\,C. Smith},  \textsc{S.\,J. Kim},  and
  \textsc{S.~Seelecke},
A stress-dependent hysteresis model for ferroelectric materials,
 \jr{Journal of Intelligent Material Systems and Structures} \textbf{18},
  69--88 (2007).


\bibitem{bassiouny89:2}
 \textsc{E.~Bassiouny} and  \textsc{A.\,F. {Ghaleb}},
Thermodynamical formulation for coupled electromechanical hysteresis effects:
  Combined electromechanical loading,
 \jr{International Journal of Engineering Science} \textbf{27}(8), 989--1000
  (1989).


\bibitem{Belov06}
 \textsc{A.\,Y. {Belov}} and  \textsc{W.\,S. {Kreher}},
Simulation of microstructure evolution in polycrystalline
  ferroelectrics-ferroelastics,
 \jr{Acta Materialia} \textbf{54}, 3463--3469 (2006).


\othercit
\bibitem{BrokateSprekels}
 \textsc{M.~Brokate} and  \textsc{J.~Sprekels},
Hysteresis and Phase Transitions (Springer, New York, 1996).


\bibitem{bv}
 \textsc{M.~Brokate} and  \textsc{A.~Visintin},
Properties of the {P}reisach model for hysteresis,
 \jr{J. Reine Angew. Math.} \textbf{402}, 1--40 (1989).


\bibitem{Cimaa02}
 \textsc{L.~Cima},  \textsc{E.~Laboure},  and  \textsc{P.~Muralt},
Characterization and model of ferroelectrics based on experimental {P}reisach
  density,
 \jr{Review of Scientific Instruments} \textbf{73}(10) (2002).


\bibitem{DKV13}
 \textsc{D.~Davino},  \textsc{P.~Krej\v{c}\'{i}},  and  \textsc{C.~Visone},
Fully coupled modelling of magnetomechanical hysteresis through thermodynamic
  compatibility,
 \jr{Smart Materials and Structures} \textbf{22}, 095009 (14pp) (2013).


\bibitem{Buelent05}
 \textsc{B.~Delibas},  \textsc{A.~Arockiarajan},  and  \textsc{W.~Seemann},
A nonlinear model of piezoelectric polycrystalline ceramics under quasi-static
  electromechanical loading,
 \jr{Journal of Materials Science: Materials in Electronics} \textbf{16},
  507--515 (2005).


\othercit
\bibitem{froehlich01:1}
 \textsc{A.~{Fr\"ohlich}},
{Mikromechanisches Modell zur Ermittlung effektiver Materialeigenschaften von
  piezoelektrischen Polykristallen},
Dissertation, {Universit\"at} Karlsruhe (TH), Forschungszentrum Karlsruhe,
  2001.


\bibitem{HKKL08}
 \textsc{T.~Hegewald},  \textsc{B.~Kaltenbacher},  \textsc{M.~Kaltenbacher},
  and  \textsc{R.~Lerch},
{ Efficient modeling of ferroelectric behaviour for the analysis of
  piezoceramic actuators},
 \jr{Journal of Intelligent Material Systems and Structures} \textbf{19}(10),
  1117--1129 (2008).


\bibitem{Huber06}
 \textsc{J.\,E. Huber},
Micromechanical modelling of ferroelectrics,
 \jr{Current Opinion in Solid State and Materials Science} \textbf{9}, 100--106
  (2006).


\bibitem{Huber01}
 \textsc{J.\,E. {Huber}} and  \textsc{N.\,A. {Fleck}},
Multi-axial electrical switching of a ferroelectric: theory versus experiment,
 \jr{Journal of the Mechanics and Physics of Solids} \textbf{49}, 785--811
  (2001).


\bibitem{hughes95:1}
 \textsc{D.\,C. {Hughes}} and  \textsc{J.\,T. {Wen}},
Preisach modeling and compensation for smart material hysteresis,
 \jr{Proc. SPIE, Active Materials and Smart Structures} \textbf{2427}, 50--64
  (1995).


\bibitem{KKHL10}
 \textsc{M.~Kaltenbacher},  \textsc{B.~Kaltenbacher},  \textsc{T.~Hegewald},
  and  \textsc{R.~Lerch},
{ Finite Element Formulation for Ferroelectric Hysteresis of Piezoelectric
  Materials},
 \jr{Journal of Intelligent Material Systems and Structures} \textbf{21},
  773--785 (2010).


\bibitem{kamlah01}
 \textsc{M.~Kamlah},
Ferroelectric and ferroelastic piezoceramics – modeling of electromechanical
  hysteresis phenomena,
 \jr{Continuum Mechanics and Thermodynamics} \textbf{13}(4), 219--268 (2001).


\bibitem{kamlah01:1}
 \textsc{M.~Kamlah} and  \textsc{U.~B\"ohle},
Finite element analysis of piezoceramic components taking into account
  ferroelectric hysteresis behavior,
 \jr{International Journal of Solids and Structures} \textbf{38}, 605--633
  (2001).


\othercit
\bibitem{KrasPokr}
 \textsc{M.~Krasnoselskii} and  \textsc{A.~Pokrovskii},
Systems with Hysteresis (Springer, Heidelberg, 1989).


\bibitem{KraynyukovaNesenenko13}
 \textsc{N.~Kraynyukova} and  \textsc{S.~Nesenenko},
Measure-valued solutions for models of ferroelectric material behavior,
 \jr{arXiv:1301.4071 [math.AP]} (2013).


\bibitem{mathz}
 \textsc{P.~Krej\v{c}\'{\i}},
Hysteresis and periodic solutions of semilinear and quasilinear wave equations,
 \jr{Math. Z.} \textbf{193}, 247--264 (1986).


\bibitem{max}
 \textsc{P.~Krej\v{c}\'{\i}},
On {M}axwell equations with the {P}reisach hysteresis operator: the
  one-dimensional time-periodic case,
 \jr{Apl.Mat.} \textbf{34}(5), 364--374 (1989).


\othercit
\bibitem{KrejciBuch}
 \textsc{P.~Krej\v{c}\'{i}},
{Hysteresis, Convexity, and Dissipation in Hyperbolic Equations} (Gakkotosho,
  Tokyo, 1996).


\bibitem{adk}
 \textsc{P.~Krej\v{c}\'{\i}},  \textsc{M.~Al~Janaideh},  and
  \textsc{F.~Deasy},
Inversion of hysteresis and creep operators,
 \jr{Physica B: Condensed Matter} \textbf{407}(9), 1354--1356 (2012).


\othercit
\bibitem{kuhnen01:1}
 \textsc{K.~Kuhnen},
{Inverse Steuerung piezoelektrischer Aktoren mit Hysterese-, Kriech- und
  Superpositionsoperatoren},
Dissertation, {Universit\"at} des Saarlandes, {Saarbr\"ucken}, 2001.


\bibitem{landis04:1}
 \textsc{C.\,M. Landis},
Non-linear constitutive modeling of ferroelectrics,
 \jr{Current Opinion in Solid State and Materials Science} \textbf{8}, 59--69
  (2004).


\bibitem{Linnemann09}
 \textsc{K.~Linnemann},  \textsc{S.~Klinkel},  and  \textsc{W.~Wagner},
A constitutive model for magnetostrictive and piezoelectric materials,
 \jr{International Journal of Solids and Structures} \textbf{46}, 1149--1166
  (2009).


\othercit
\bibitem{mayergoyz91:2}
 \textsc{I.\,D. Mayergoyz},
{Mathematical Models of Hysteresis} (Springer-Verlag New York, 1991).


\bibitem{McMeeking07}
 \textsc{R.\,M. McMeeking},  \textsc{C.\,M. Landis},  and  \textsc{S.\,M.\,A.
  Jimenez},
A principle of virtual work for combined electrostatic and mechanical loading
  of materials,
 \jr{International Journal of Non-Linear Mechanics} \textbf{42}(6), 831--838
  (2007).


\bibitem{MieheZaehRosato12}
 \textsc{C.~Miehe},  \textsc{D.~Z\"ah},  and  \textsc{D.~Rosato},
Variational-based modeling of micro-electro-elasticity with electric
  field-driven and stress-driven domain evolutions,
 \jr{International Journal for Numerical Methods in Engineering} \textbf{91},
  115--141 (2012).


\bibitem{MielkeTimofte06}
 \textsc{A.~Mielke} and  \textsc{A.~Timofte},
An energetic material model for time-dependent ferroelectric behavior:
  Existence and uniqueness,
 \jr{Math. Methods Appl. Sci.} \textbf{29}, 1393--1410 (2006).


\bibitem{Pasco04}
 \textsc{Y.~Pasco} and  \textsc{A.~Berry},
A hybrid analytical/numerical model of piezoelectric stack actuators using a
  macroscopic nonlinear theory of ferroelectricity and a {P}reisach model of
  hysteresis,
 \jr{Journal of Intelligent Material Systems and Structures} \textbf{15},
  375--386 (2004).


\bibitem{preisach}
 \textsc{F.~Preisach},
{\"U}ber die magnetische {N}achwirkung,
 \jr{Z. Physik} \textbf{94}, 277--302 (1935).


\othercit
\bibitem{SchroederKeip10}
 \textsc{J.~Schr\"oder} and  \textsc{M.\,A. Keip},
Multiscale modeling of electro--mechanically coupled materials: homogenization
  procedure and computation of overall moduli,
 in: Proceedings of the IUTAM conference on multiscale modeling of fatigue,
  damage and fracture in smart materials,  (Springer, Heidelberg, 2010).


\bibitem{Schroeder05}
 \textsc{J.~Schr\"oder} and  \textsc{H.~Romanowski},
A thermodynamically consistent mesoscopic model for transversely isotropic
  ferroelectric ceramics in a coordinate-invariant setting,
 \jr{Archive of Applied Mechanics} \textbf{74}, 863--877 (2005).


\bibitem{Smith12}
 \textsc{R.\,C. Smith} and  \textsc{Z.~Hu},
The homogenized energy model (hem) for characterizing polarization and strains
  in hysteretic ferroelectric materials: Material properties and uniaxial model
  development,
 \jr{Journal of Intelligent Material Systems and Structures} \textbf{23},
  1833--1867 (2012).


\bibitem{Smith03}
 \textsc{R.\,C. Smith},  \textsc{S.~Seelecke},  \textsc{Z.~Ounaies},  and
  \textsc{J.~Smith},
A free energy model for hysteresis in ferroelectric materials,
 \jr{Journal of Intelligent Material Systems and Structures} \textbf{14},
  719--737 (2003).


\bibitem{Landis07}
 \textsc{Y.~Su} and  \textsc{C.\,M. Landis},
Continuum thermodynamics of ferroelectric domain evolution: Theory, fnite
  element implementation and application to domain wall pinning,
 \jr{Journal of the Mechanics and Physics of Solids} \textbf{55}, 280--305
  (2007).


\othercit
\bibitem{VisintinBuch}
 \textsc{A.~Visintin},
Differential Models of Hysteresis (Springer, Berlin, 1994).


\bibitem{Wangetal}
 \textsc{J.~Wang},  \textsc{M.~Kamlah},  and  \textsc{T.\,Y. Zhang},
{Phase field simulations of low dimensional ferroelectrics},
 \jr{Acta Mechanica} \textbf{214}, 49--59 (2010).


\bibitem{Xuetal}
 \textsc{B.\,X. Xu},  \textsc{D.~Schrade},  \textsc{R.~M\"uller},
  \textsc{D.~Gross},  \textsc{T.~Granzow},  and  \textsc{J.~R\"odel},
{Phase field simulation and experimental investigation of the
  electro-mechanical behavior of ferroelectrics},
 \jr{Z. Angew. Math. Mech.} \textbf{90}, 623--632 (2010).


\end{thebibliography}
%

\end{document}